\documentclass[11pt]{article}
\usepackage{authblk} 
 \usepackage{tikz}
 \usetikzlibrary{shapes}

\usepackage{lipsum}
\usepackage{amsfonts}
\usepackage{amsmath,amsopn}
\usepackage{graphicx}
\usepackage{epstopdf}
\usepackage{algorithm}
\usepackage{algpseudocode}
\ifpdf
  \DeclareGraphicsExtensions{.eps,.pdf,.png,.jpg}
\else
  \DeclareGraphicsExtensions{.eps}
\fi

\usepackage{amssymb}
\usepackage{tabularx, booktabs}
\usepackage{tabularx}

\newtheorem{example}{Example}[section]
\usepackage{subcaption}

\DeclareMathOperator{\MSE}{MSE} 
\DeclareMathOperator{\Swish}{Swish}

\usepackage[colorlinks=true, allcolors=blue]{hyperref}

\usepackage[top=2.8cm,bottom=2.8cm,left=2.5cm,right=2.5cm]{geometry}
\usepackage{float}

\usepackage[font=footnotesize,width=\linewidth]{caption} 

\newtheorem{theorem}{Theorem}
\newtheorem{lemma}{Lemma}

\title{Alikhanov--XfPINNs: Adaptive Physics-Informed Learning for Nonlinear Fractional PDEs on Nonuniform Meshes}

\author{Himanshu Kumar Dwivedi\footnote{ \href{himanshukrdwivedi.rs.mat21@itbhu.ac.in}{himanshukrdwivedi.rs.mat21@itbhu.ac.in}} ,
Matthias Ehrhardt\footnote{Corresponding author,  \href{mailto:ehrhardt@uni-wuppertal.de}{ehrhardt@uni-wuppertal.de}} ,
Rajeev\footnote{ \href{rajeev.apm@iitbhu.ac.in}{rajeev.apm@iitbhu.ac.in}} 
}

\affil{Department of Mathematical Sciences, IIT (BHU), Varanasi 221005, India}
\affil{IMACM, School of Mathematics and Natural Sciences, \\ University of Wuppertal, Germany}

\begin{document}
\maketitle

\begin{tikzpicture}[remember picture,overlay]
	\node[anchor=north east,inner sep=20pt] at (current page.north east)
	{\includegraphics[scale=0.2]{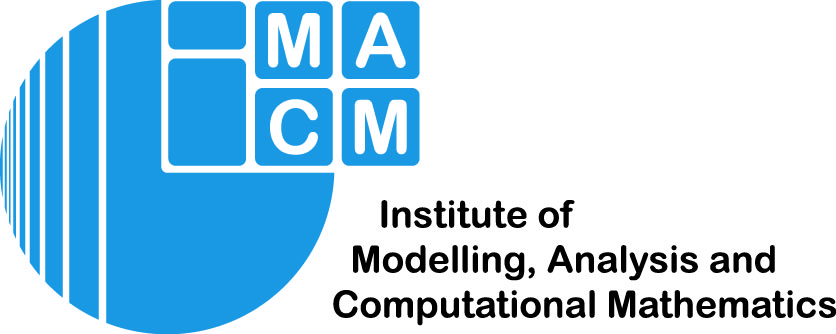}};
\end{tikzpicture}

\begin{abstract}
To address the initial singularity inherent in solutions to fractional partial differential equations (fPDEs), we propose an accelerated Alikhanov discretization formulation implemented on nonuniform time grids. 
Based on the \textit{physics-informed neural networks (PINNs)} framework, 
we introduce an Alikhanov-extended fractional PINNs (XfPINNs) architecture that combines high-order temporal discretization and deep learning. 
The nonlocal memory term in fPDEs leads to high computational cost, while the weak singularity near $t\to 0^+$ can deteriorate accuracy on uniform meshes. 
To separate temporal discretization effects from optimization and sampling errors, we further develop an auxiliary time-marching configuration that enables auditable temporal-convergence studies under controlled training tolerances.
This architecture can solve general nonlinear fPDEs. 
The XfPINNs approach is designed for forward and inverse problems, allowing for data-driven solution reconstruction and parameter estimation. 
First, the neural network approximates the solution of nonlinear fPDEs; then, an adaptive activation function accelerates convergence and enhances training efficiency. 
The optimization framework embeds a variational loss function constructed from the Alikhanov scheme, 
where the initial and boundary conditions are imposed using a combination of hard and soft constraints. 
Numerical experiments, including cases with known and unknown exact solutions which demonstrate the robustness, computational efficiency, and significant CPU time savings of the Alikhanov--XfPINNs method.
\end{abstract}

\begin{minipage}{0.9\linewidth}

\medskip

\end{minipage}

%%%%%%%%%%%%%%%%%%%%%%%%%%%%%%%% Intro
\section{Introduction} % ME: rewritten
The integration of classical scientific computing and modern machine learning has given rise to the emerging field of \textit{scientific machine learning}, offering a paradigm shift for approximating partial differential equations (PDEs). 
Leveraging abundant data and advances in computational power, recent machine learning developments have driven remarkable progress across disciplines, from image recognition \cite{P11} to cognitive science \cite{P12}. 
At the forefront of this advancement is the \textit{physics-informed neural network} (PINN) framework \cite{P21,P22}, which harnesses the universal approximation capabilities of deep neural networks to solve fractional PDEs by embedding the governing physical laws directly into the learning architecture.

PINNs are appealing because they are versatile and easy to implement, making them highly effective for solving diverse classes of PDEs.
Their applications span computational physics, fluid dynamics \cite{P27, P26}, meta-material design \cite{P213}, and biomedical modeling \cite{P24}. 
Due to their mesh-free nature and reliance on randomly sampled collocation points, PINNs are well-suited for high-dimensional problems  \cite{P215}
and can handle complex geometries and irregular domains with ease \cite{P219,P220}.  
In particular, \textit{fractional PINNs} (fPINNs) have emerged as a promising approach for modeling anomalous diffusion and memory-driven phenomena prevalent in fractional-order systems. 
By incorporating discretizations of fractional-order derivatives, such as the Riemann-Liouville or Caputo operators,
into the loss function, fPINNs can accurately solve and identify parameters in nonlocal and weakly singular PDEs. 
Recent advancements include combining fPINNs with graded temporal schemes and fast convolution approximations 
to overcome the challenges posed by long-memory effects and singular behavior near the initial time.

In this work, we consider a nonlinear, parametrized time-fractional PDE of the general form
\begin{equation}\label{caputo}
{} _0^{\mathcal{C}}\mathcal{\partial}_t^{\alpha} v+\mathcal{N}[v; \lambda]=0,
\quad t\in [0,T], \;\mathbf{x}\in \Omega, 
\end{equation}
where $\partial_t^{\alpha}=\;_{0}^C \mathcal{D}_t^{\alpha}$ denotes the \textit{Caputo time-fractional derivative of order $\alpha$} 
 \begin{equation}
{} _0^{\mathcal{C}}\mathcal{\partial}_t^{\alpha} u(t):=\bigl((\mathcal{RL})_t^{(1-\alpha)} (u)'\bigr)(t) 
=\int_0^t \omega_{1-\alpha}(t-\xi)\, (u(\xi))^{'}\,d\xi, 
\quad 0<\alpha<1,
\end{equation}
where $(\mathcal{RL})_t^{\beta}$ is defined below and referred as Riemann-Liouville fractional integral of order $\beta>0$,
\begin{equation}
    (\mathcal{RL})_t^{\beta}:= \int_0^t \omega_{\beta}(t-\xi)u(\xi)\,d\xi,
    \quad \omega_{\beta}(s)=\frac{s^{\beta-1}}{\Gamma(\beta)}. 
\end{equation}
% where $\omega_{\beta}(s)=\frac{s^{\beta-1}}{\Gamma(\beta)}$. 
The function $v(\cdot, t)$ represents the latent solution, and $\mathcal{N}[\cdot; \lambda]$ denotes a nonlinear differential operator parameterized by the vector $\lambda$. 
This formulation encompasses a wide range of mathematical models in physics, including conservation laws, diffusion processes, advection-reaction-diffusion systems, and kinetic equations. 
For instance, the % one-dimensional 
1D time-dependent viscous Burgers equation \cite{P113} corresponds to the case $\mathcal{N}[v; \lambda] = \lambda_1 v v_x - \lambda_2 v_{xx}$, where $\lambda = (\lambda_1, \lambda_2)$.  Classical numerical approaches such as finite-difference, finite-element, and spectral methods \cite{P410} -- have long been employed to address forward and inverse problems for \emph{nonlinear fractional partial differential equations} (NfPDEs).

In this work, we develop a second-order Alikhanov scheme on nonuniform temporal grids, guided by the following key considerations
\begin{itemize}
\item Time-fractional PDEs with Caputo derivatives exhibit a weak singularity near $t=0$, cf.\ \cite{jin,stynes2,lean}. 
Thus, graded meshes are recommended for accurately resolving the initial layer \cite{stynes2,stynes1}.
\item 

Classical schemes for time-fractional PDEs typically require $\mathcal{O}(MK_t)$ memory and $\mathcal{O}(MK_t^2)$ operations, where $M$ and $K_t$ denote the spatial and temporal grid sizes, respectively. To reduce this burden, we develop an accelerated Alikhanov scheme inspired by \cite{sjiang2017,dwivedi1,dwivedi_2,fast_tempered,secondfast,fast_non}.
\end{itemize}
The goal of inverse problems in fPDEs is to identify unknown functions or parameters using limited observational data, thereby revealing key mechanisms in physical, chemical, or biological systems. 
Traditional methods, such as regularization, eigenfunction expansions, and Laplace transforms \cite{P417,P418}, are widely used, but they have limitations in high-dimensional settings and are sensitive to complex initial and boundary conditions. 
Additionally, they lack algorithmic flexibility.
Inverse problems are commonly classified into two types: function reconstruction \cite{P419} and parameter estimation \cite{P420}. 
The latter is the focus of this work in the context of NfPDEs.

PINNs \cite{P21} provide a powerful alternative by embedding the governing physical laws directly into the neural network loss function during training. 
This allows PINNs to solve both forward and inverse fPDEs \cite{P425,P427} with enhanced generalizability, 
high-dimensional scalability \cite{P429}, and robustness to noisy or sparse data \cite{P4}.
However, a key challenge lies in the nonlocal nature of fractional derivatives, such as Caputo and Riemann–Liouville derivatives, 
which are not directly amenable to automatic differentiation. 
Several fPINN variants have been proposed to overcome this challenge by incorporating finite difference \cite{P431},
spectral, and power series \cite{P433} discretizations into the loss function. 
These models enable accurate solution recovery and parameter inference. 
Further developments include hybrid interpolation schemes \cite{P435}, dual-network structures \cite{P434}, and applications to complex fPDEs with oscillatory or singular dynamics.

Recent progress in fPINNs has focused on incorporating numerical schemes for Caputo derivatives into neural networks. 
However, challenges remain, particularly regarding the high cost due to memory effects, initial singularities,
and suboptimal training efficiency in large-scale problems. 
To address these issues, we present an accelerated Alikhanov-fPINNs approach on nonuniform time grids 
for high-dimensional time-fractional PDEs.
The key contributions are:
\begin{itemize}
\item To improve training efficiency, we impose initial and boundary conditions through suitable hard and soft constraint formulations, covering both homogeneous and nonhomogeneous cases. We further incorporate an adaptive activation function with a trainable parameter $(a)$ and a scaling factor $(n)$.

\item The proposed method is validated across diverse benchmark tests, 
including high-dimensional forward and inverse problems, showcasing its robustness and computational effectiveness.
\end{itemize}

%%% the structure
The manuscript is organised as follows: 
Section~\ref{sec2} introduces the neural network framework for nonlinear fractional PDEs. 
Section~\ref{xfpinn_framew} derives the Alikhanov approximation and its accelerated version on nonuniform time grids. 
Section~\ref{sect_aap} presents the implementation of Alikhanov-based fPINNs/XfPINNs with adaptive activation functions. 
Section~\ref{sec5} reports numerical results for forward and inverse problems in multi-dimensional settings. 
The final conclusions are given in Section~\ref{sec6}.

%%%%%%%%%%%%%%%%%%%%%%%%%%%%%%%%%%%%%%%%%%%% SECTION 2 %%%%%%%%%%%%
\section{Neural Network Framework for Fractional PDEs}\label{sec2}
First, we revisit the core PINN methodology as applied to NfPDEs. 
Then, in Section~\ref{xfpinn_framew}, we present our Alikhanov--XfPINN framework for nonuniform time meshes of time-fractional NfPDEs.
We provide a detailed global consistency error analysis for the accelerated sum-of-exponentials (SOE) approximation on nonuniform meshes, explicitly accounting for the initial singularity.

%%%%%%%%%%%%%%% Section 2.1.
\subsection{Physics-Informed Neural Networks}
PINNs provide an effective framework for approximating complex PDEs, including models with nonlocal or fractional dynamics. Their main strength lies in embedding the governing physical laws directly into the training objective, thereby producing neural approximations that remain consistent with the underlying equations rather than relying solely on data. In practice, PINNs can be implemented efficiently using standard deep learning libraries such as TensorFlow and PyTorch.

First, consider a general nonlinear fractional PDE (NfPDE) defined over a spatial domain $\Omega\subset\mathbb{R}^d$, a time interval $t\in(0,T]$, and a fractional order $\alpha \in (0,1)$. 
The problem is formulated through the governing operator $(\mathcal{G})$, boundary operator $(\mathcal{B})$, and initial operator $(\mathcal{I})$ as
\begin{align*}
\mathcal{G}[\mathbf{x}, t, v(\mathbf{x}, t; \lambda)] &= 0,
\quad (\mathbf{x},t) \in \Omega \times (0,T],\\
\mathcal{B}[\mathbf{x}, t, v(\mathbf{x}, t; \lambda)]& = 0,
\quad \mathbf{x} \in \partial\Omega,; t\in(0,T],\\
\mathcal{I}[\mathbf{x}, v(\mathbf{x},0; \lambda)]& = 0,
\quad \mathbf{x} \in \Omega .
\end{align*}

Here, $v(\mathbf{x},t)$ denotes the solution, and the operator $\mathcal{G}$ typically involves a Caputo-type fractional derivative $\partial_t^\alpha v$ along with nonlinear spatial components, as $\nabla v, \nabla^2 v,\dots$.

Mathematically, a deep neural network can be viewed as a nested composition of multivariate functions. 
A PINN seeks to approximate the mapping $v_{NN}(\mathbf{x},t)\colon\Omega \times [0,T]\to\mathbb{R}$, which associates each input pair $(\mathbf{x},t)\in\Omega\times[0,T]$ with the scalar output $v_{NN}(\mathbf{x},t)\in\mathbb{R}$ through a sequence of linear operations and nonlinear activations.  
The resulting representation for $1\le l\le L-1$ is given by
\begin{equation}\label{eq:NN_layers}
\begin{aligned}
\text{input layer :} \quad &NN^0 = (\mathbf{x}, t),\\
\text{hidden layers :} \quad &NN^l = \sigma\bigl(w^l \cdot NN^{l-1} + b^l\bigr),\quad \\
\text{output layer :} \quad &NN^L = w^L \cdot NN^{L-1} + b^L.
\end{aligned}
\end{equation}
Here, $NN^l$ denotes the output of the $l$-th layer.
$w$ and $b$ are the trainable weight matrices and bias vectors, respectively, and $\sigma$ is the activation function.  
For each hidden layer ($l\ge1$), the previous layer's output is first mapped by a learnable affine transformation and then passed through $\sigma$, and finally forwarded to the next layer. 
The final layer only applies the affine mapping to produce the network’s output. 
This architecture is purely feed-forward, with no recurrent or feedback connections.

% In physics‐informed neural networks,
In PINNs, the governing PDE is incorporated explicitly into the loss functional, ensuring that the learned solution satisfies the underlying physics. To illustrate this framework, we consider a general form of a NfPDE.
\begin{equation}\label{ac_eq}
\begin{split}
{} _0^{\mathcal{C}}\mathcal{\partial}_t^{\alpha} v
+\mathcal{N}[v;\lambda] &= g(\mathbf{x},t),
\quad t\in [0,T],\, x\in\Omega,\\
v(\mathbf{x},0) &= \varphi(\mathbf{x}),
\qquad\mathbf{x}\in\Omega,\\
v(\mathbf{x},t) &= \phi(\mathbf{x},t),\quad(\mathbf{x},t)\in\partial\Omega\times[0,T],
\end{split}
\end{equation}
where $v(\mathbf{x}a,t)$ denotes the exact analytical solution, $\varphi(\mathbf{x})$ specifies the initial conditions, $\phi(\mathbf{x},t)$ enforces the boundary constraints and $\partial\Omega$ marks the spatial boundary of $\Omega$. 
The PINNs embeds the NfPDE~\eqref{ac_eq} by minimizing the residual 
between the PDE’s left‐ and right‐hand sides. 
All required derivatives are obtained via automatic differentiation to ensure a consistent application of the chain rule.
Let the training data be denoted by
\begin{equation*}
   \mathcal{I}={\mathcal{I}_f,\,\mathcal{I}_{ic},\,\mathcal{I}_{bc}},
\end{equation*}
where $(\mathcal{I}_f)$, $(\mathcal{I}_{ic})$, and $(\mathcal{I}_{bc})$ represent the interior collocation, initial-condition, and boundary-condition datasets, respectively:
\begin{align*}
\mathcal{I}_f&=\bigl\{(\mathbf{x}_f^i,t_f^i,g(\mathbf{x}_f^i,t_f^i))\bigr\}_{i=1}^{N_f},\\
\mathcal{I}_{ic}&=\bigl\{(\mathbf{x}_{ic}^i,0,\varphi(\mathbf{x}_{ic}^i))\bigr\}_{i=1}^{N_{ic}},\\
\mathcal{I}_{bc}&=\bigl\{(\mathbf{x}_{bc}^i,t_{bc}^i,\phi(\mathbf{x}_{bc}^i,t_{bc}^i))\bigr\}_{i=1}^{N_{bc}}.
\end{align*}
Accordingly, the total loss associated with PINN is defined as
\begin{equation}\label{pinn_loss}
     \mathcal{L}(\Theta) = \MSE_f+\MSE_{ic}+\MSE_{bc}.
\end{equation}
where the mean-squared error terms corresponding to the PDE residual, initial condition, and boundary constraints are defined as follows:
\begin{align*}
   \MSE_{f} &= \frac{1}{N_{f}}\sum_{i=1}^{N_f}
\biggl(\frac{\partial v_{NN}(\mathbf{x}_f^i, t_f^i)}{\partial t} +\mathcal{N}[v_{NN}(\mathbf{x}_f^i, t_f^i);\lambda ]
- g(\mathbf{x}_f^i, t_f^i)\biggr)^2,\\
\MSE_{ic} &= \frac{1}{N_{ic}}\sum_{i=1}^{N_{ic}}
\bigl(v_{NN}(\mathbf{x}_{ic}^i,0)-\varphi(\mathbf{x}_{ic}^i)\bigr)^2,\\
\MSE_{bc} &= \frac{1}{N_{bc}}\sum_{i=1}^{N_{bc}}
\bigl(v_{NN}(\mathbf{x}_{bc}^i, t_{bc}^i) -\phi(\mathbf{x}_{bc}^i, t_{bc}^i)\bigr)^2.
\end{align*}
%In the above equations, $\Theta$ denotes the unknown parameters $w$ and $b$ of the networks. $v_{NN}$ represents the output of neural network. The loss function measures the extent to which the PINN solution satsifies the Allen-Cahn equations and IBCs. Obviously, there is an error between the approximated solution $v_{NN}$ and the real solution $v(\varkappa,t)$, but the error can be reduced by adjusting the parameters $\Theta$. Therefore, the goal of this algorithm is to find a set of parameters $\Theta$ such that the approximated solution $v_{NN}$ minimizes $\mathcal{L}(\Theta)$. If it is relatively small, then the approximated solution $v_{NN}$ is considered to fit equation~\eqref{ac_eq}. This process of obtaining optimal parameters of the neural network is called "training". The training process enables the neural network to construct a good mapping between input and output for Allen-Cahn equation.
Here, $\Theta$ denotes the complete set of trainable network parameters, including weights and biases. The objective functional $\mathcal{L}(\Theta)$ quantifies how well the neural approximation $v_{\mathrm{NN}}$ satisfies the NfPDE together with the prescribed initial and boundary data. Although $v_{\mathrm{NN}}$ is generally not identical to the exact solution $(v(\mathbf{x},t))$, the mismatch can be reduced by optimizing $\Theta$. Thus, the training process seeks
\begin{equation*}
  \Theta^*=\arg\min_{\Theta}\mathcal{L}(\Theta),
\end{equation*}
hence, $v_{\mathrm{NN}}$ produces a sufficiently small residual in \eqref{ac_eq} and satisfies the imposed constraints with high accuracy. 
This optimization procedure constitutes the training stage, during which the network learns an effective approximation to the NfPDE solution.

%%%%%%%%%%%%%%%%%%%%%%%%%%%%%%% Section 2.2. %%%%%%%%%%%%%%%%
\subsection{Extension to Alikhanov--XfPINNs on Nonuniform Time Meshes}\label{xfpinn_framew}
Motivated by the strengths of PINNs, we propose Alikhanov--XfPINNs, a generalized framework designed to solve nonlinear fractional PDEs on nonuniform time grids. 
A central challenge in solving such problems is the nonlocal nature of the Caputo derivative, which does not permit the direct application of the chain rule. 
To overcome this challenge, we use a second-order Alikhanov scheme to discretize the temporal fractional derivative on a nonuniform mesh and integrate it directly into the loss function.

To accelerate convergence and enhance representation power, we also incorporate adaptive activation functions. 
Although this study focuses on the Caputo derivative, the Alikhanov--XfPINNs framework can be adapted to other types of fractional operators (e.g., the Riemann–Liouville derivative) by modifying the corresponding discretization technique. 
This flexibility makes the proposed methodology a versatile tool for solving a broad class of fractional PDEs.

%%%%%%%%%%%%%%%%%%%%%%%%%%%%%%%%%%%%% Section 2.2.1
 \subsubsection{\textit{Time Discretization on Nonuniform Grids}}\label{graded_structure}
Throughout this work, we employ nonuniform temporal discretizations. Let
\begin{equation*}
0=t_0<t_1<\cdots<t_{k-1}<t_k<\cdots<t_{K_t}=T
\end{equation*}
be a partition of the time interval $[0,T]$, with local time step
\begin{equation*}
\tau_k:=t_k-t_{k-1},\qquad 1\le k\le K_t,
\end{equation*}
and maximum step size
%\begin{equation*}
$\tau:=\max_{1\le k\le K_t}\tau_k$.
%\end{equation*}
For $\theta=\alpha/2\in(0,1)$, the off-set time level is defined as
\begin{equation*}
t_{k-\theta}:=\theta t_{k-1}+(1-\theta)t_k.
\end{equation*}
The adjacent step ratio and its global bound are denoted by
\begin{equation*}
\rho_k:=\frac{\tau_k}{\tau_{k+1}},\qquad
\rho:=\max_{k\ge1}\rho_k.
\end{equation*}

We impose the following mesh conditions.
%%%%%%%%%%%%%%%% M1
\begin{description}\label{M1}
\item[$\mathbf{M1}$:] The maximum step ratio satisfies
\begin{equation*}
   \rho\le \frac{3}{2}.
\end{equation*}
\end{description}

Assumption $\mathbf{M1}$ permits moderate variation between consecutive time steps. In particular, it allows time steps to decrease by at most a factor of $2/3$, while no additional restriction is imposed on their growth. Such a condition is standard in variable-step fractional discretizations, although an upper bound on $\tau$ may still be required in the stability analysis.

Nonuniform meshes are especially useful for fractional evolution problems, since their solutions often show weak singular behavior around $t=0$, for example
\begin{equation*}
u_t=\mathcal{O}(t^{\alpha-1})\qquad \text{as }t\to0,
\end{equation*}
even for smooth data. Graded meshes of the form
\begin{equation*}
t_n=T\Big(\frac{n}{K_t}\Big)^\gamma
\end{equation*}
are therefore commonly used to cluster time levels near $t=0$ and improve accuracy; see, for example, \cite{stynes2}. Additional mesh-regularity assumptions are needed to ensure convergence on general nonuniform grids \cite{liao16}.
%%%%%%%%%%%%%%%%%%%% M2
\begin{description}\label{M2}
\item[$\mathbf{M2}$:] There exist constants $C_{1\gamma},C_{2\gamma}>0$, independent of the mesh size and a grading parameter $\gamma\ge1$ such that
\begin{equation*}
\tau_n\le \tau\min\{1,C_{1\gamma}t_n^{1-1/\gamma}\},
\qquad 1\le n\le K_t,
\end{equation*}
and
\begin{equation*}
t_n\le C_{2\gamma}t_{n-1},
\qquad 2\le n\le K_t.
\end{equation*}
\end{description}

The qunatity $\gamma$ determines the degree of refinement near the origin. The choice $\gamma=1$ corresponds to a quasi-uniform mesh satisfying $\mathbf{M2}$, while larger values of $\gamma$ produce stronger clustering near $t=0$ and are better suited for resolving initial singularities.

%%%%%%%%%%%%%%%%%%%%%%%%%%%%%% Section 2.2.2.
\subsubsection{\textit{Alikhanov Formula for a Nonuniform Time Mesh}}
We consider a grid function $\{v^k\}_{k=0}^K$ on a nonuniform time mesh. 
Next, we define the weighted operator 
$v^{k-\theta}=\theta v^{k-1}+(1-\theta)v^k$ 
and the standard difference operator $\triangledown_{\tau} v^k=v^k-v^{k-1}$.
To streamline the presentation and clarify the notation, we introduce:
\begin{equation*}
\bar{\omega}_k=\omega_{2-\alpha}(t-t_{k-\theta}), 
\quad\text{and}\quad 
    \bar{\omega}_k'=\omega_{1-\alpha}(t_{k-\theta}-t),\quad 
    t\in[0,t_{k-\theta}].
\end{equation*} 
Let $(\mathcal{L}_{1,k}v)$ denote the linear interpolant of $v$ at the nodes $(t_{k-1})$ and $(t_k$), while $(\mathcal{Q}_{2,k}v)$ represents the quadratic interpolant constructed from $(t_{k-1}), (t_k)$, and $(t_{k+1})$.

%\end{itemize}
We define
\begin{equation}
   (\mathcal{L}_{1,k}v)'(t)=\frac{\triangledown_{\tau} v^k }{\tau_k}, 
   \quad(\mathcal{Q}_{2,k}v)'(t)=\frac{\triangledown_{\tau}v^k}{\tau_k}+\frac{2(t-t_{k-1/2})}{\tau_k(\tau_k+\tau_{k+1})}(\rho_k\triangledown_{\tau}v^{k+1}-\triangledown_{\tau}v^k).
\end{equation}
We now present the \textit{nonuniform Alikhanov formula}, which reads:
\begin{equation}\label{eq3.1_P4}
\begin{split}
({} _0^{\mathcal{C}}\mathcal{\partial}_t^{\alpha}v)^{k-\theta}&=\int_{t_{k-1}}^{t_{k-\theta}} \omega_{1-\alpha}(t_{k-\theta}-s)\,(\mathcal{L}_{1,k}v)'(s)\,ds\\
& \qquad+\sum_{n=1}^{k-1} \int_{t_{n-1}}^{t_n} \omega_{1-\alpha}(t_{k-\theta}-s) \,(\mathcal{Q}_{2,n}v)'(s)\,ds\\
&= \boldsymbol{\mathrm{a}}^{(0,k)} \triangledown_{\tau} v^k+\sum_{n=1}^{k-1} \bigl(\boldsymbol{\mathrm{a}}^{(k-n,k)}\triangledown_{\tau}v^n+\boldsymbol{\mathrm{b}}^{(k-n,k)} (\rho_n \triangledown_{\tau} v^{n+1}-\triangledown_{\tau}v^n)\bigr).
\end{split}
\end{equation}
We define
\begin{align}
  \boldsymbol{\mathrm{a}}^{(k-n,k)} &= \frac{1}{\tau_n}\int_{t_n}^{\min\{t_n,t_{k-\theta}\}} \bar{\omega}_k'(s) \,ds,
  \quad 1\le n\le k,\\
  \boldsymbol{\mathrm{b}}^{(k-n,k)} &= 2\int_{t_{n-1}}^{t_n} \frac{(s-t_{n-1/2}) }{\tau_n(\tau_n+\tau_{n+1})}\bar{\omega}_k'(s)\,ds,\quad  1\le n\le k-1.
\end{align}
Rearranging the terms in \eqref{eq3.1_P4}, we obtain the compact form 
\begin{equation}
     ({} _0^{\mathcal{C}}\mathcal{\partial}_t^{\alpha}v)^{k-\theta}=\sum_{n=1}^k \mathcal{D}^{(k-n,k)}\triangledown_{\tau} v^n,
\end{equation}
where the discrete convolution kernels are given by $\mathcal{D}^{(0,1)}=a^{(0,1)}$ if $k=1$, and 
\begin{align}\label{eq3.5_P4} 
   \mathcal{D}^{(k-n,k)}  =\left\{
    \begin {aligned}
         &  \boldsymbol{\mathrm{a}}^{(k-1,k)} -\boldsymbol{\mathrm{b}}^{(k-1,k)},   \quad & n=1,  \\
         & \boldsymbol{\mathrm{a}}^{(k-n,k)} +\rho_{n-1}\boldsymbol{\mathrm{b}}^{(k-n+1,k)} -\boldsymbol{\mathrm{b}}^{(k-n,k)},  \quad &2\le n\le k-1,\\
           &   \boldsymbol{\mathrm{a}}^{(0,k)} +\rho_{k-1}\boldsymbol{\mathrm{b}}^{(1,k)},  \quad & n=k  .\\              
    \end{aligned}
\right.
\end{align}
The nonuniform formula \eqref{eq3.1_P4} generalizes the Alikhanov scheme on uniform meshes \cite{Alikhanov} and the nonuniform counterpart was first introduced in \cite{fast_non} to address initial singularities via graded meshes. 
The following recent developments concerning the properties of the coefficients of accelerated approximation can be found in \cite{Liao_arxiv}. 
%%%%%%%%%%%%%%%%%%%%%%%%%%%%%%%%%%%% Lemma 2.1.
\begin{lemma}[\cite{Liao_arxiv}]\label{main_3.1}
Assume that the time-step condition $\mathbf{M1}$ is satisfied. We then examine the discrete convolution kernels $\mathcal{D}^{(k-n,k)}$ defined in \eqref{eq3.5_P4}. The following properties hold:

\begin{enumerate}
\item[\textbf{(a)}] For $1\le n\le k$, the kernels $\mathcal{D}^{(k-n,k)} $ satisfy \label{main_a}
\begin{equation}
\begin{split}
 \frac{24}{11\tau_{k}} \int_{t_{k-1}}^{t_k} \omega_{1-\alpha}(t_k-s)\,ds &\ge  \mathcal{D}^{(0,k)},\\
  \mathcal{D}^{(k-n-1,k)}-\mathcal{D}^{(k-n,k)} &
  \ge\frac{4}{11\tau_{k}}\int_{t_{k-1}}^{t_k}\omega_{1-\alpha}(t_n-s)\,ds. 
  \end{split}
\end{equation}
%%%%%%%
 \item[\textbf{(b)}] When $1\le n\le k-1$, the kernels $\mathcal{D}^{(k-n,k)}$ are monotone \label{main_b} 
 \begin{equation}\label{monotone}
\begin{split}
  0&\le (1+\rho_k)\boldsymbol{\mathrm{b}}^{(k-n,k)}-\frac{1}{5\tau_n}\int_{t_{n-1}}^{t_n} (t_n-s)\omega_{1-\alpha}(t_{k-\theta}-s)\,ds\\
  &<\mathcal{D}^{(k-n-1,k)}-\mathcal{D}^{(k-n,k)} ,
  \end{split}
\end{equation}
%%%%%%%%%%%%%
  \item[\textbf{(c)}] For $k\ge 2$, the leading kernel $\mathcal{D}^{(0,k)}$ dominates the subsequent kernel, namely,\label{main_c}
\begin{equation}
     \mathcal{D}^{(1,k)}< \frac{(1-2\theta)}{(1-\theta)}\,\mathcal{D}^{(0,k)}.
  \end{equation}
\end{enumerate}
\end{lemma}

%%%%%%%%%%%%%%%%%%%%%%%%%%%%%%%% Section 2.2.3
\subsubsection{\textit{Accelerated Alikhanov formula for nonuniform time mesh}}
The approximation \eqref{eq3.5_P4} becomes expensive in long-time simulations because the fractional derivative retains the full solution history. To reduce this memory and computational burden, we use the sum-of-exponentials approximation \cite{sjiang2017} for the kernel, yielding an accelerated approximation. The main idea is as follows.

%%%%%%%%%%%%%%%%
\begin{lemma}[\cite{sjiang2017}]\label{SOE_lemma} 
When $0<\alpha<1$, the absolute tolerance error $\epsilon \ll 1$, a cut-off time size $\Delta t >0$, and a final time $T$ are given, then there exists a positive integer $\mathcal{N}_q$, a set of positive quadrature nodes $s^l$, and their respective positive weights $\nu^l(1\le l\le\mathcal{N}_q)$, with
\begin{equation}
    \Bigl| \omega_{1-\alpha}(t)-\sum_{l=1}^{\mathcal{N}_q} \nu^l\,\mathrm{e}^{-s^l t}\Bigr|
    \le \epsilon,
    \quad \forall t\in[\Delta t, T],
\end{equation}
and the number $\mathcal{N}_q$ of required quadrature nodes is of order 
\begin{equation}
    \mathcal{N}_q =
    O\biggl(
    \log\frac{1}{\epsilon}
    \Bigl(\log\log\frac{1}{\epsilon}+\log\frac{T}{\Delta t}\Bigr)
    +
    \log\frac{1}{\Delta t}
    \Bigl(\log\log\frac{1}{\epsilon}+\log\frac{1}{\Delta t}\Bigr)
    \biggr).
\end{equation}
\end{lemma}

As shown in Lemma~\ref{SOE_lemma}, we decompose the Caputo derivative \eqref{caputo} into a local term over $[t_{k-1},t_{k-\theta}]$ and a history term over $[0, t_{k-1}]$.
The local term is approximated using linear interpolation of $v’(t)$.
The history term is efficiently evaluated via the sum-of-exponentials approximation of $\bar{\omega}_k’(s)$, yielding:
\begin{equation}\label{eq5.1_P4}
\begin{split}
({} _{0}^{\mathcal{C}}\mathcal{\partial}_{t}^{\alpha}v)^{k-\theta}
&\approx \int_{t_{k-1}}^{t_{k-\theta}}\bar{\omega}_k'(s) (\mathcal{L}_{1,k}v)'(s)\,ds
+\int_{0}^{t_{k-1}}\sum_{l=1}^{\mathcal{N}_q}\nu^l\,\mathrm{e}^{-s^l(t_{k-\theta}-s)}v'(s)\,ds,\\
&=\boldsymbol{a}^{(0,k)}\triangledown_{\tau}v^k+\sum_{l=1}^{\mathcal{N}_q}\nu^l\mathcal{V}_{his}^l(t_{k-1}),\quad k\ge1,
\end{split}
\end{equation}
where $\mathcal{V}_{his}^l(t_0)=0$ 
and $\mathcal{V}^l_{his}(t_k) = \int_0^{t_k} \mathrm{e}^{-s^l(t_{k-\theta}-\xi)}v'(\xi)\,d\xi$.
To evaluate $(\mathcal{V}_{\mathrm{his}}^l(t_k))$, we employ a recursive update obtained by approximating $v'(t)$ with a quadratic interpolant on each subinterval $[t_{l-1},t_l]$, $1\le l\le k-1$, as follows:

\begin{align}
    \mathcal{V}^l_{his}(t_n)
    &\approx \int_0^{t_{n-1}} \mathrm{e}^{-s^l(t_{n+1-\theta}-s) }v'(s)\,ds
    +\int_{t_{n-1}}^{t_n} \mathrm{e}^{-s^l(t_{n+1-\theta}-s)}(\mathcal{Q}_{2,n}v)'(s)\,ds,\nonumber\\
&= \mathrm{e}^{-s^l(\theta \tau_n+(1-\theta) \tau_{n-1})}\mathcal{V}^l_{his}(t_{n-1})+\boldsymbol{c}^{(n,l)}\triangledown_{\tau}v^n +\boldsymbol{d}^{(n,l)} (\rho_n\triangledown_{\tau}v^{n+1} -\triangledown_{\tau}v^n),\label{eq5.2_P4}
\end{align}
where the discrete coefficients are defined by
\begin{align}
\boldsymbol{c}^{(n,l)} &= \frac{1}{\tau_n}\int_{t_{n-1}}^{t_n}
\mathrm{e}^{-s^l(t_{n+1-\theta} -\xi) }\,d\xi,\\
\boldsymbol{d}^{(n,l)} &= \int_{t_{n-1}}^{t_n} \mathrm{e}^{-s^l(t_{n+1-\theta}-\xi )}\frac{2(\xi-t_{n-1/2}) }{\tau_n(\tau_n+\tau_{n+1})}\,d\xi.
\end{align}
From \eqref{eq5.1_P4}--\eqref{eq5.2_P4}, we obtain the following form of the approximation
\begin{equation}\label{fast_eq}
({} _{0}^{\mathcal{FC}}\mathcal{\partial}_{t}^{\alpha} v)^{k-\theta}=\boldsymbol{a}^{(0,k)}\triangledown_{\tau}v^k+\sum_{l=1}^{\mathcal{N}_q}\nu^l\mathcal{V}^l_{his}(t_{k-1}), \quad k\ge1.
\end{equation}
Here, $\mathcal{V}^{l}_{\mathrm{his}}(t_k)$ is evaluated through the following recurrence relation.

\begin{equation}\label{recursive_1}
\mathcal{V}^l_{his}(t_n) = \mathrm{e}^{-s^l(\theta \tau_n+(1-\theta)\tau_{n+1})}\mathcal{V}^l_{his}(t_{n-1})+\boldsymbol{c}^{(n,l)}\triangledown_{\tau} v^n+\boldsymbol{d}^{(n,l)}(\rho_n \triangledown_{\tau}v^{n+1} -\triangledown_{\tau}v^n).
\end{equation}
By Lemma~\ref{SOE_lemma}, the number of quadrature nodes satisfies $\mathcal{N}_q=\mathcal{O}(\log N)$ for $T\ge 1$ and $\mathcal{N}_q=\mathcal{O}(\log^2 N)$ for $T\approx 1 $, leading to substantial reductions in computational cost and memory consumption.

We examine the scheme by defining the consistency error associated with the accelerated approximation \eqref{eq5.1_P4}:
\begin{equation*}
\Upsilon^j[u]
:=({}_{0}^{\mathcal{C}}\partial_t^\alpha u)(t_{j-\theta})
-({}_{0}^{\mathcal{FC}}\partial_t^\alpha u)^{j-\theta},
\qquad j\ge1.
\end{equation*}
Following the framework of \cite{Liao_arxiv}, Lemma~\ref{lemma23} establishes a discrete convolution-type bound for $\Upsilon^j[u]$ on general nonuniform time meshes.
Furthermore, the fractional Gr\"onwall lemma shows that the dominant contribution to the error in solution come from the convolution part $\sum_{j=1}^n \mathsf{C}^{(n-j,n)} |\Upsilon^j[u]|$, here  complementary kernels $\mathsf{C}^{(k-n,k)}$ defined as 
 \begin{equation}\label{p_des2}
 \begin{split}
  \mathsf{C}^{(0,k)} &:= \frac{1}{\mathcal{D}^{(0,k)}},\\
  \mathsf{C}^{(k-n,k)} &:=\frac{1}{\mathcal{D}^{(0,n)}} \sum_{i=n+1}^{k}(\mathcal{D}^{(i-n-1,i)} -\mathcal{D}^{(i-n,i)} )\mathsf{C}^{(k-i,k)},\quad 1\le n\le k-1.
\end{split}
\end{equation}
% Moreover, Lemma \ref{lem3.1_P4}’s fractional Grönwall inequality indicates that the solution error is primarily driven by the convolution term $\sum_{j=1}^n C(n)_{n-j} |\Upsilon^j[u]|$, where $C(n)_{n-j}$ are the complementary kernels defined in \eqref{p_des}.
%Drawing inspiration from \cite{graded_liao}, we establish in the subsequent lemma that $\Upsilon_{j-\theta}$ can be effectively governed via a discrete convolution framework, applicable to a wide spectrum of temporal meshes. Furthermore, the fractional Grönwall inequality presented in Lemma \ref{lem3.1_P4} demonstrates that the error in the solution is primarily influenced by the convolution error $\sum_{j=1}^n C(n)_{n-j} |\Upsilon^j[u]|$, where $C(n)_{n-j}$ represent the complementary kernels described in \eqref{p_des}.

%%%%%%%%%%%%%%%%%%%%%%%% Lemma 2.3.
\begin{lemma}[\cite{fast_non}]\label{lemma23} 
Suppose that the condition $\mathbf{M1}$ on the step ratio is satisfied, and that $v \in \mathcal{C}^3(0,T]$ with $\int_0^T \xi^2 |v'''(\xi)|\,d\xi < \infty$, 
moreover a constant $\mathcal{C}_v>0$ such that
\begin{equation*}
  |v'''(t)|\le \mathcal{C}_v\bigl(1+t^{\upsilon-2}\bigr),\qquad 0<t\le T,
\end{equation*}
where $\upsilon\in(0,1)$ characterizes the temporal regularity. For the fast Alikhanov formula on nonuniform grids \eqref{eq5.1_P4} with the discrete convolution kernels, denoted as $\mathcal{D}^{(k-n,k)}$, along with the local consistency error $\Upsilon^{j-\theta}$, a convolutional framework is exhibited, 
as illustrated below
\begin{multline}
  |\Upsilon^{k-\theta}[v]|
  \le\mathcal{D}^{(0,k)}\mathcal{G}_{loc}^k 
  +\sum_{i=1}^{k-1} \Bigl(\mathcal{D}^{(k-i-1,k)}-\mathcal{D}^{(k-i,k)} \Bigr)\mathcal{G}_{his}^i
+\frac{\mathcal{C}_v}{\upsilon}\hat{t}_{k-1}^2\epsilon,\quad 1\le k\le K_t,
\end{multline}
where the two quantities $\mathcal{G}_{loc}$ and $\mathcal{G}_{his}$ are
\begin{align*}
  \mathcal{G}_{loc}^k &= \frac{3}{2}\biggl(\int_{t_{k-1}}^{t_{k-1/2}} (s-t_{k-1})^2 |v'''(s)|\,ds 
+\tau_k \int_{t_{k-1/2}}^{t_{k-1}}(t_k-s)|v'''(s)|\,ds   \biggr),\\
\mathcal{G}_{his}^k &= \frac{5}{2}\biggl(\int_{t_{k-1}}^{t_{k}} (s-t_{k-1})^2 |v'''(s)|\,ds
+ \int_{t_{k}}^{t_{k+1}}(t_{k+1}-s)^2|v'''(s)|\,ds\biggr),
\end{align*}
and $\hat{t}_k = \max\limits_{1\le k\le {K_t}}\{1,t_k\}$. 
Then, the global convolution error is 
\begin{equation}
  \sum_{i=1}^{k} \mathsf{C}^{(k-i,k)} |\Upsilon^i[v]|
   \le\sum_{i=1}^{k}\mathsf{C}^{(k-i,k)}   \mathcal{D}^{(0,k)}\mathcal{G}_{loc}^{k}+\sum_{i=1}^{k} \mathsf{C}^{(k-i,k)} \mathcal{D}^{(0,k)}\mathcal{G}_{his}^{k}+\frac{\mathcal{C}_v}{\upsilon}t_{k}^{\alpha}\hat{t}_{k-1}^2\epsilon.
\end{equation}
\end{lemma}
%% the following better at the beginning
% In conclusion, we present an in-depth global consistency analysis of errors for the approximation \eqref{fast_eq} when applied to nonuniform meshes, taking into account the initial singularity.
Lemma~\ref{lemma23} directly yields the following theorem, which establishes a global error bound for the accelerated approximation.

%%%%%%%%%%%%%%%%%%%%%%%% Theorem 2.4
\begin{theorem}[\cite{Liao_arxiv}]\label{err_thm}
Suppose that $v\in \mathcal{C}^3(0,T]$ and a constant $\mathcal{C}_v>0$ such that
\begin{equation*}
|v'''(t)|\le \mathcal{C}_v\bigl(1+t^{\upsilon-2}\bigr),
\quad 0<t\le T,
\end{equation*}
where temporal regularity $\upsilon\in(0,1)$. 
Under the assumptions stated in $\mathbf{M1}$,
the global consistency error of the accelerated approximation \eqref{fast_eq} satisfies the following estimate:

\begin{multline}
\sum_{i=1}^{k}\mathsf{C}^{(k-i,k)} |\Upsilon^i[v]|
\le \mathcal{C}_v\bigg(\frac{\tau_1^{\upsilon}}{\upsilon}+\frac{1}{(1-\alpha)}\max\limits_{2\le n\le k}t_{n}^{\alpha}t_{n-1}^{\upsilon-3}\tau_{n}^3\tau_{n-1}^{-\alpha}+\frac{\epsilon}{\upsilon}t_{k}^{\alpha}\hat{t}_{k-1}^2\bigg),\quad 1\leq n\leq K_t.
\end{multline}
Particularly, under a temporal discretization satisfying the graded-type condition $\mathbf{M2}$, it follows that
\begin{equation}
   \sum_{i=1}^{k}\mathsf{C}^{(k-i,k)}  |\Upsilon^i[v]|\leq  \frac{\mathcal{C}_v}{\upsilon(1-\alpha)}\bigg( \tau^{\min (\gamma \upsilon,2)}+\epsilon \bigg).
\end{equation}
\end{theorem}
%The subsequent lemma from \cite{graded_liao} is provided to 
%The subsequent theorem % from \cite{Liao_arxiv}  
%reveals that error in temporal direction from the term $v^{k-\theta}$ is constrained by  Alikhanov approximation error.
The following theorem shows that the error in the temporal direction of the term $v^{k-\theta}$ is limited by the Alikhanov approximation error.

%%%%%%%%%%%%%%%%%%%%%%% Theorem 2.5.
\begin{theorem}[\cite{Liao_arxiv}]\label{err_thm2}
Let $v\in C^2((0,T])$ and a positive constant $\mathcal{C}_v$ such that
\begin{equation*}
    |v''(t)|\le \mathcal{C}_v\bigl(1+t^{\upsilon-2}\bigr),
    \qquad 0<t\le T,
\end{equation*}
where $\upsilon\in(0,1)$ denotes the temporal regularity parameter. The local truncation error associated with $v^{k-\theta}$ is defined by
\begin{equation*}
    \mathcal{R}^k[v]:=v(t_{k-\theta})-v^{k-\theta},\qquad 1\le k\le K_t.
\end{equation*}
Under the graded-type mesh condition $\mathbf{M2}$, the corresponding global consistency error is bounded for $1\leq k\leq K_t$ as follows:
\begin{equation}
   \sum_{i=1}^k \mathsf{C}^{(k-i,k)}  |\mathcal{R}_{\tau}^i[v]|\le
\mathcal{C}_v\biggl(\frac{\tau_1^{\upsilon+\alpha}}{\upsilon} +t_k^{\alpha} \max\limits_{2\le n\le k}t_{n-1}^{\upsilon-2}\tau_{n}^2  \biggr)
\end{equation}
is taken into account. 
\end{theorem}

%Inspired by physics-informed neural networks (PINNs), we introduce the $\mathcal{F{L}}_{21}\sigma$-XfPINN methodology on graded meshes for solving TFAC equation. Since the Caputo fractional derivative does not admit a direct application of the chain rule, we employ an accelerated $\mathcal{F{L}}2-1\sigma$ approximation on a graded mesh to discretize this operator and incorporate it into the loss function. To further enhance training efficiency, we integrate an adaptive activation function together with a residual-based, stage-wise adaptive training strategy. Although this work focuses on the Caputo derivative, the $\mathcal{F{L}}_{21}\sigma$-XfPINN framework can be extended to other fractional derivatives—such as the Riemann–Liouville derivative—by suitably adapting the discretization scheme. This flexibility broadens the applicability of $\mathcal{F{L}}_{21}\sigma$-XfPINN to a wide class of fractional partial differential equations (FPDEs).

%%%%%%%%%%%%%%%%%%%%%%%%%%%%%%%% Section 2.3. %%%%%%
\subsection{Adaptive Activation Function}\label{sect_aap}
\begin{figure}[htb]
  \centering
  \makebox[\textwidth][c]{%
  \includegraphics[width=12.9cm, height=8.75cm]{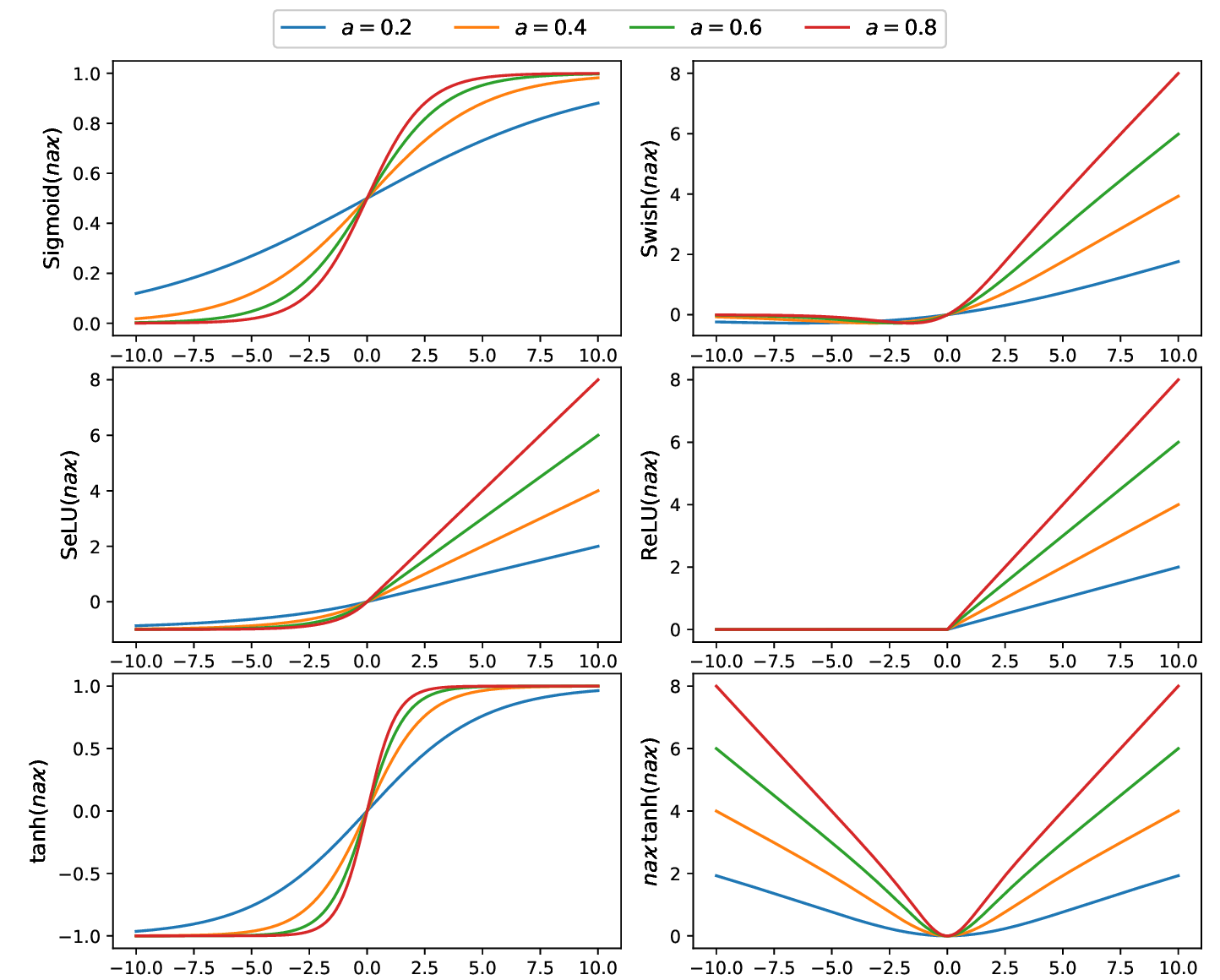}%
  }
  \caption{Standard activation functions modified adaptively by varying the hyperparameter $a$, with the scaling factor fixed at $n=1$.}
  \label{fig_six_adap_f}
\end{figure}

To optimize the Alikhanov--XfPINN on graded time grids more quickly, 
each standard activation is replaced by its adaptive counterpart. 
his counterpart is parameterized by a trainable slope $a$ and a prescribed scaling factor $n$, as proposed by Jagtap, Kawaguchi, and Karniadakis~\cite{jagtap}.
The six adaptive activation functions considered in this work are listed in Table~\ref{tab:activation}. 
The parameter $a$ adjusts the local steepness of the activation and thus modifies the gradient flow during optimization, as illustrated in Figure ~\ref{fig_six_adap_f}. 
The additional scaling factor $n$ amplifies this adaptive effect and helps mitigate slow convergence caused by small learning rates.

For the Alikhanov--XfPINN framework on a nonuniform temporal mesh, the neural approximation is written as
\begin{equation}
v_{NN}(\mathbf{x},t;\Theta)
=\bigl(\phi_h\circ\sigma\circ n a\,\phi_{h-1}
\circ\cdots\circ\sigma\circ n a\,\phi_{1}\bigr)(\mathbf{x},t),
\end{equation}
where $\circ$ denotes function composition and
\begin{equation*}
   \phi_i = W_i\phi_{i-1}+b_i,\qquad i=1,\ldots,h.
\end{equation*}
Here, the spatio-temporal input is $\phi_0=(\mathbf{x},t)$, 
and the resulting network output is $v_{NN}(\mathbf{x},t;\Theta)$.
In the forward $\mathcal{FL}_{21}\sigma$-XfPINN setting, the trainable parameter set is
\begin{equation*}
   \Theta = \{W_i,b_i,a\}_{i=0}^h,
\end{equation*}
whereas, in the inverse Alikhanov--XfPINN formulation, the unknown physical parameter $\zeta$ is also learned, giving
\begin{equation*}
   \Theta=\{W_i,b_i,a,\zeta\}_{i=0}^h.
\end{equation*}

%%%%%%%%%%%%%%%%%%%%%%%%%%%%% Table 1
\begin{table}[ht]
  \centering
  \renewcommand{\arraystretch}{2.2}
  \caption{Performance comparison: standard vs.\ adaptive variants of activation functions.}
  \label{tab:activation}
  \begin{tabularx}{\textwidth}{@{}l *{2}{>{\centering\arraybackslash}X}@{}}
    \toprule
    \textbf{Activation} & \textbf{Fixed} & \textbf{Adaptive} \\
    \midrule
    Sigmoid
      & $\displaystyle \frac{1}{1 + \mathrm{e}^{-\mathbf{x}}}$
      & $\displaystyle \frac{1}{1 + \mathrm{e}^{-na\,\mathbf{x}}}$ \\

    Swish
      & $\displaystyle \frac{\mathbf{x}}{1 + \mathrm{e}^{-\mathbf{x}}}$
      & $\displaystyle \frac{na\,\mathbf{x}}{1 + \mathrm{e}^{-na\,\mathbf{x}}}$ \\

    SeLU
      & $\displaystyle
        \begin{cases}
          \lambda\,\mathbf{x}, & \mathbf{x} > 0,\\[6pt]
          \lambda\,\alpha\,(\mathrm{e}^{\mathbf{x}}-1), & \mathbf{x}\le 0
        \end{cases}$
      & $\displaystyle
        \begin{cases}
          \lambda\,na\,\mathbf{x}, & \mathbf{x} > 0,\\[6pt]
          \lambda\,\alpha\,(\mathrm{e}^{na\,\mathbf{x}}-1), & \mathbf{x} \le 0
        \end{cases}$ \\

    ReLU
      & $\max(0,\mathbf{x})$
      & $\max(0,\,na\,\mathbf{x})$ \\

    tanh
      & $\tanh(\mathbf{x})$
      & $\tanh(na\,\mathbf{x})$ \\

    $\mathbf{x}$\,tanh
      & $\mathbf{x}\,\tanh(\mathbf{x})$
      & $na\,\mathbf{x}\,\tanh(na\,\mathbf{x})$ \\
    \bottomrule
  \end{tabularx}
\end{table}

%%%%%%%%%%%%%%%%%%%%%%%%%%%%%% Figure 1

%%%%%%%%%%%%%%%%%%%%%%%%%%%%%% Section 2.4.
\subsection{Loss Function}
We construct the proposed framework by embedding the prescribed constraints directly into the trial solution. 
Unlike penalty-based approaches to initial and boundary conditions, this strategy incorporates the prescribed constraints directly into the neural ansatz \cite{geonet}. 
Consequently, separate penalty terms for these conditions are not required in the loss functional, and fewer training samples are needed for constraint enforcement.

\begin{figure}[htb]
\hspace{0.09\textwidth} % Adjust this value to fine-tune %horizontal position
%\begin{subfigure}{0.8\textwidth}
   \includegraphics[width=0.8\textwidth]{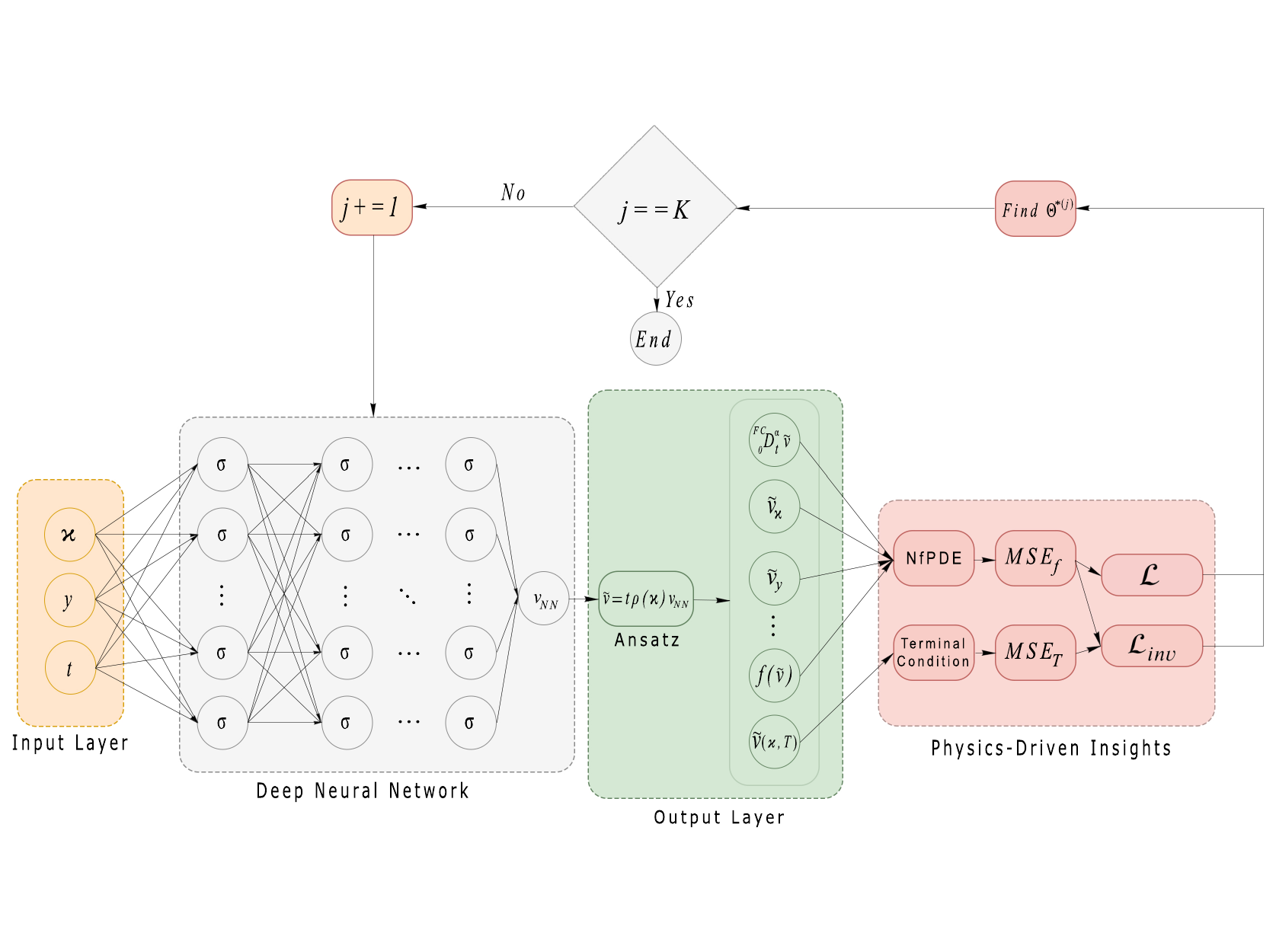}
%  \caption{}
%  \label{gph51}
%\end{subfigure}
\caption{Schematic diagram of the XfPINN architecture.}
\label{gph:test5}
\end{figure}

Specifically, we define the constrained neural approximation by
\begin{equation}\label{ansatz}
  \tilde{v}(\mathbf{x},t)
  =  t\,\rho(\mathbf{x})\,v_{NN}(\mathbf{x},t;\Theta),
\end{equation}
where $v_{NN}(\mathbf{x},t;\Theta)$ denotes the neural-network output and $\rho(\mathbf{x})$ is a given boundary mask chosen to satisfy
\begin{equation*}
   \rho(\mathbf{x})=0,\qquad \mathbf{x}\in\partial\Omega.
\end{equation*}
The factor $t$ in \eqref{ansatz} enforces the initial condition, while $\rho(\mathbf{x})$ imposes the homogeneous boundary constraint.
Therefore, the approximation $\tilde{v}(\mathbf{x},t)$ satisfies the prescribed initial and boundary conditions (IBCs) by construction. 
Figure~\ref{gph:test5} depicts the proposed Alikhanov–XfPINN structure under a graded time discretization, where the 
% initial and boundary conditions 
IBCs are enforced through the hard-constraint ansatz.

Since the IBCs
% initial and boundary conditions 
are already embedded in the ansatz, the training set only contains interior residual points
\begin{equation*}
    \mathcal{I} = \{\mathcal{I}_f\},\qquad
     \mathcal{I}_f = \{(\mathbf{x}_f^i,t_f^i,g(\mathbf{x}_f^i,t_f^i))\}_{i=1}^{N_f}.
\end{equation*}
Here, $N_f$ represents the total number of residual sampling points.
Therefore, the forward-training objective only contains the NfPDE residual term:
\begin{equation}\label{xfpinn_loss}
    \mathcal{L}(\Theta) = \MSE_f,
\end{equation}
with
\begin{equation*}
  \MSE_f = \frac{1}{N_f}\sum_{i=1}^{N_f}
  \Bigl[{} _{0}^{\mathcal{C}}\mathcal{\partial}_t^{\alpha} \tilde{v}(\mathbf{x}_f^i,t_f^i) + \mathcal{N}(\tilde{v}(\mathbf{x}_f^i,t_f^i);\lambda)- g(\mathbf{x}_f^i,\,t_f^i) \Bigr]^2.
\end{equation*}
In the inverse Alikhanov--XfPINN formulation, we augment the forward‐loss with a terminal‐time term so that the NfPDEs parameters and network weights are learned simultaneously. 
We define the terminal dataset  
$\mathcal{I}_T = \bigl\{(\mathbf{x}_T^i,\,t_T^i,\,\psi(\mathbf{x}_T^i))\bigr\}_{i=1}^{N_T}$,  
where $N_T$ is its size. 
With the combined dataset $\mathcal{I} = \{\mathcal{I}_f,\,\mathcal{I}_T\}$, the hard‐constraint inverse loss is  
\begin{equation}\label{xfpinn_inv_loss}
     \mathcal{L}_{\rm inv}(\Theta)
     = w_f\MSE_f +\, w_T\MSE_T,
\end{equation}
where
\begin{equation*}
  \MSE_T = \frac{1}{N_T}
  \sum_{i=1}^{N_T}\Bigl[\tilde{v}(\mathbf{x}_T^i,T) - \psi(\mathbf{x}_T^i)\Bigr]^2.
\end{equation*}
The learnable parameters are $\Theta = \{W_i,\,b_i,\,a,\,\zeta\}_{i=0}^h$,  
and $w_f$, $w_T$ are the respective weights for $\MSE_f$ and $\MSE_T$.

%%%%%%%%%%%%%%%%%%%%%%%%%%%%%%%%%% Section 2.5.
\subsection{Adaptive Training Process}
Pang, Lu, and Karniadakis~\cite{P431} developed fractional PINNs for space-time fractional advection-diffusion models. They treated fractional derivatives using finite-difference discretizations and computed integer-order differential terms via automatic differentiation.  
For models involving Caputo time-fractional derivatives, they employed the standard $\mathcal{L}1$ discretization and determined the trainable parameters $\Theta^*$ by optimizing the physics-informed loss.

To reduce the cost associated with the nonlocal memory of the Caputo operator and better resolve the weak singularity near $t=0$, Shi, Liu, and Yang~\cite{shi} developed the FL1--fPINN method on graded meshes. 
This approach combines an accelerated nonuniform $\mathcal{L}1$ approximation with the fPINN framework to efficiently simulate multi-dimensional time-fractional PDEs.

In this work, we focus on models governed by time-fractional dynamics. 
While fPINNs and FL1--fPINNs can both be applied to the time-fractional Allen–Cahn (TFAC) problem, our numerical results demonstrate that the Alikhanov--XfPINN framework offers higher convergence rates and significantly greater accuracy across the benchmark tests considered.

First, we replace the accelerated $\mathcal{L}1$ formula with an accelerated Alikhanov approximation on a graded mesh to improve the convergence rate. 
As Theorems~\ref{err_thm} and \ref{err_thm2} establish, this approximation achieves order~2 for smooth solutions.
For solutions with an initial singularity, the convergence rate is $\min\{\gamma\alpha,2\}$, where $\gamma$ denotes the grading parameter. 
Depending on the value of $\gamma$, the convergence order can be as high as~2. 
Thus, by minimizing the loss function, we can significantly improve performance.

Second, we introduce a residual-based progressive training strategy that improves computational efficiency and stabilizes the optimization process. 
At stage $j$, the residual training set is restricted to the current temporal window and is defined as
\begin{equation*}
   \mathcal{I}_f^j = \bigl\{
   (\mathbf{x}_f^i,t_f^i,g(\mathbf{x}_f^i,t_f^i))\in\mathcal{I}_f
   \colon t_f^i\le t_{j-\theta}\bigr\},
   \qquad j=1,\ldots,K-1 .
\end{equation*}
The training process begins with the initial short-time interval $[t_{1-\theta}, t_{2-\theta}]$, 
in which the network uses the dataset $I_f^1$ to learn the local approximation $\tilde{v}^1$. 
For the subsequent stages ($j=3,\ldots,K$), 
the temporal domain is gradually expanded to the interval
$[t_{1-\theta}, t_{j-\theta}]$, and the network is retrained with an updated dataset $\mathcal{I}_f^j$.
In the final stage, the network uses the dataset $\mathcal{I}_f^K$ to approximate the function $\tilde{v}$
over the entire computational interval, $[0, t_{K-\theta}]$.
Each stage is stopped once the prescribed tolerance is reached or the maximum number of iterations is attained. 
At stage $j$, the trainable parameters are obtained by solving
\begin{equation}\label{loss_XfPINN}
\mathcal{L}^j(\Theta)
= \frac{1}{N_f^j}
  \sum_{(\mathbf{x}_f,t_f,g)\in\mathcal{I}_f^j}
\Bigl[{} _{0}^{\mathcal{C}}\mathcal{\partial}_t^{\alpha}\tilde{v}(\mathbf{x}_f,t_f) + \mathcal{N}(\tilde{v}(\mathbf{x}_f,t_f);\lambda)- g\Bigr]^2,
\end{equation}
yielding $\Theta_j^* = \arg\min\mathcal{L}^j(\Theta)$, where $N_f^j = |\mathcal{I}_f^j|$. 
%To further enhance the stage-wise training efficiency, we embed a Residual-based Adaptive Refinement (RAR) [cite deepxde] procedure within each training stage. After an initial set of $N_{\rm init}$ optimizer iterations on the stage-$j$ training set $\mathcal{I}_f^j$, we perform $N_{\rm RAR}$ refinement cycles. In each cycle, we first form a \emph{candidate pool} by uniformly sampling (with replacement) $P$ points from the current stage set $\mathcal{I}_f^j$, ensuring all times satisfy $t_i\le t_{j+\sigma}$.  Here $P=\max(c\,|\mathcal{I}_f^j|,\,N_{\rm new})$ with candidate factor $c$ and minimal new-sample count $N_{\rm new}$. We then evaluate the absolute PDE residual on the pool and select the top $N_{\rm new}$ points with largest residuals. These hardest points are appended to the stage training set, and training continues for another $N_{\rm upd}$ iterations to drive down both data and PDE losses. This RAR cycle is repeated $N_{\rm RAR}$ times before advancing to the next time-stage. By adaptively mining the regions where the network residual is largest, RAR focuses computational effort on the most challenging parts of the solution, leading to accelerated convergence and improved accuracy with fewer total training epochs.
%
We follow the exact same stage-wise schedule as in the inverse Alikhanov--XfPINN formulation, 
but now we minimize the minimize the inverse loss with a hard-constraint at stage $j$:
\begin{equation}\label{loss_XfPINN_inv}
\begin{split}
\mathcal{L}_{\rm inv}^j(\Theta)&= w_f\,\frac{1}{N_f^j}\sum_{(\mathbf{x}_f,t_f,g)\in\mathcal{I}_f^j} 
\Bigl[{} _{0}^{\mathcal{C}}\mathcal{\partial}_t^{\alpha} \tilde{v}(\mathbf{x}_f,t_f) + \mathcal{N}(\tilde{v}(\mathbf{x}_f,t_f);\lambda)- g \Bigr]^2\\
&\qquad + w_T\,\frac{1}{N_T}  \sum_{(\mathbf{x}_T^i,t_T^i,\psi)\in\mathcal{I}_T}
\Bigl[\tilde{v}(\mathbf{x}_T^i,T)-\psi(\mathbf{x}_T^i)\Bigr]^2.
\end{split}
\end{equation}
Because $\mathcal{I}_f^j$ only contains points with $t_f\le t_{j-\theta}$, the second (terminal) term vanishes for all $j<K$ and only “turns on” at the final stage $j=K$.
All subsequent steps, including the implementation of residual-based adaptive refinement on $\mathcal{I}_f^j$ are performed identically to the forward Alikhanov--XfPINN formulation. 

The proposed framework is implemented in TensorFlow~\cite{tensorflow} to efficiently evaluate the residual terms and optimize the trainable parameters. 
TensorFlow provides automatic differentiation to compute spatial derivatives and gradient information, optimized CPU/GPU kernels to efficiently train models, direct GPU support, and a flexible Python interface to rapidly implement and test models. 
The broader TensorFlow ecosystem, which includes visualization and monitoring tools such as TensorBoard, facilitates systematic model development.

We employ the Adam algorithm~\cite{adam} for optimization, 
with suitable hyperparameter adjustments for the graded-mesh Alikhanov--XfPINN setting.
These adjustments are made to improve convergence speed while maintaining numerical accuracy. 
The complete training strategy is summarized below.
Algorithm~\ref{alg:xfpinn} summarizes the computational workflow of the Alikhanov--XfPINN for NfPDEs. 
First, we introduce the temporal nodes $t_k$, where $k = 0,\dots,K$.
Then, we construct the initial dataset $\mathcal{I}$ by randomly selecting spatio-temporal samples in $\Omega_T$. 
In the stage-wise Alikhanov--XfPINN procedure, the active training set is updated at each stage according to the current temporal window (see line~\ref{line_I} of Algorithm~\ref{alg:xfpinn}. 
Once the collocation samples are fixed, we use the accelerated approximation to compute the fractional derivative 
and evaluate the total loss. Then, we update the trainable parameters using gradient-based optimization.

In Algorithm~\ref{alg:xfpinn}, the trainable variables are only initialized during the first stage $(r=1)$. 
For all subsequent stages, training begins with the optimal parameters obtained in the previous stage. 
This implementation uses TensorFlow's automatic differentiation in Python 
and adopts Xavier initialization for the network weights to ensure stable signal propagation.
The neural output is embedded into the hard-constrained form \eqref{ansatz}, producing the complete approximation $v_{\Theta}(x,y,t)$. 
At each collocation point corresponding to $t_k$, the network values at the preceding time levels $\{t_j\}_{j=0}^k$
are evaluated simultaneously in vectorized form. This allows the fractional history term to be computed efficiently.

%%%%%%%%%%%%%%%%%%%

%%%%%%%%%%%%%%%%%%%%%%%%%%
\begin{algorithm}
\scriptsize
\caption{Alikhanov--XfPINNs on nonuniform meshes for simulation and estimation in NfPDEs}
\label{alg:xfpinn}
\begin{algorithmic}[1]
  \Require A set of training dataset $\mathcal{I}$, domain bounds $\mathbf{lb},\mathbf{ub}$, network layer sizes $\{d_0,d_1,\dots,d_L\}$, fractional order $\alpha$, subdivisions $K$, shift parameter $\theta$, $\lambda_1$, $\lambda_2$, stage-wise max.\ iterations $m_{\mathrm{stage}}$ and threshold $\varepsilon$.
  \Ensure Trained parameters $\Theta^*$ and approximate solution $\tilde{v}(\mathbf{x},t)$.

  \State Initialize weights and biases:
  \For{$l = 0$ to $L-1$}
    \State Sample weight matrix $W^{(l)} \sim \mathcal{N}(0,\sigma_l^2)$ with $\sigma_l = \sqrt{2/(d_l + d_{l+1})}$ 
    \State Set bias vector $b^{(l)} \gets \mathbf{0}$
  \EndFor

  \State $\Theta^{*(0)} \gets$ initial parameters
  \For{$j = 2$ to $K$}
    \State $t_{j-\theta} \gets t_{j-1} + (1-\theta)\,\tau_{j}$
    \State\label{line_I} $\mathcal{I}_f^j \gets \{(\mathbf{x}^i,t_f^i,g(\mathbf{x}^i,t^i)) \in \mathcal{I} :t_{1-\theta} \le t_f^i \le t_{j-\theta} \}$
    \State $N_f^j\gets|\mathcal{I}_f^j|$
    \State $\Theta^{(j)} \gets \Theta^{*(j-1)}$
    \For{$k = 1$ to $m_{\mathrm{stage}}$}
      \State Compute $\tilde{v}(\mathbf{x},t)
  = t\,\rho(\mathbf{x})\,v_{NN}(\mathbf{x},t;\Theta^{(j)})$ and ${} _0^{\mathcal{C}}\mathcal{\partial}_t^{\alpha}\tilde{v}(\mathbf{x},t)$ via Alikhanov formula
      \State Compute stage loss for the forward and inverse problems:
      \begin{equation*}
        \mathcal{L}^j(\Theta^{(j)})
        = \frac{1}{N_f^j}
        \sum_{(\mathbf{x}_f,t_f,g)\in\mathcal{I}_f^j}
\sum_{i=1}^{N_f}
\Biggl[\partial_t^\alpha \tilde{v}(\mathbf{x}_f,t_f) + \mathcal{N}(\tilde{v}(\mathbf{x}_f,t_f);\lambda)- g(\mathbf{x}_f,\,t_f) \Biggr]^2
\end{equation*}
      
\begin{align*}
\mathcal{L}_{\rm inv}^j(\Theta) = w_f\,\frac{1}{N_f^j}&\sum_{(\mathbf{x}_f,t_f,g)\in\mathcal{I}_f^j} \bigl[{} _{0}^{\mathcal{C}}\mathcal{\partial}_t^{\alpha} \tilde{v}(\mathbf{x}_f,t_f) + \mathcal{N}(\tilde{v}(\mathbf{x}_f,t_f);\lambda)- g(\mathbf{x}_f,\,t_f)\bigr]^2\\
& +w_T\,\frac{1}{N_T}  \sum_{(\mathbf{x}_T^i,t_T^i,\psi)\in\mathcal{I}_T}\bigl[\tilde{v}(\mathbf{x}_T^i,T)-\psi(\mathbf{x}_T^i)\bigr]^2.
\end{align*}
      \State Update:
      \begin{equation*}
        \Theta^{(j)} \leftarrow \Theta^{(j)} - \eta\,\nabla_{\Theta^{(j)}}\,\mathcal{L}^{(j)}(\Theta^{(j)})
        \quad\text{(via {tf.keras.optimizers.Adam})}
       \end{equation*}
      \If{$\mathcal{L}^{(j)}(\Theta^{(j)}) < \varepsilon$}
        \State \textbf{break}
      \EndIf
    \EndFor
    \State Store $\Theta^{*(j)} \gets \Theta^{(j)}$
  \EndFor

  \State \Return $\Theta^* = \Theta^{*(K-1)}$, $\tilde{v}(\mathbf{x},t)
  = t\,\rho(\mathbf{x})\,v_{NN}(\mathbf{x},t;\Theta^*)$
\end{algorithmic}
\end{algorithm}
%% this is the algorithm related to Alikhanov-fPINN-M: HKD
\begin{algorithm}
\scriptsize
\caption{Alikhanov-fPINN-M on nonuniform meshes for time-marching enforcement}
\label{alg:asa_fpinn_m}
\begin{algorithmic}[1]
\Require Collocation points $\{\boldsymbol{x}_f^i\}_{i=1}^{N_x}\subset\Omega$, final time $T$, number of time levels $K$, grading parameter $\gamma_g$, offset parameter $\theta$, model data $\lambda_1$, $\lambda_2$, $\lambda_3$, $g$, physical parameters $(\alpha,\varepsilon,\gamma)$, SOE nodes and weights $\{(s^l,\nu^l)\}_{l=1}^{N_{\exp}}$, discretization coefficients $\{\boldsymbol{\mathrm{a}}^{(0,k)},a^{(n,l)},b^{(n,l)}\}$, maximum inner iterations $m_{\mathrm{step}}$, stopping tolerance $\varepsilon_{\mathrm{tol}}$
\Ensure Frozen snapshots $\{\tilde{v}^{k}\}_{k=1}^{K}$ and stored history states $\{\mathcal{V}^l(t_k;\cdot)\}_{k=1}^{K}$

\State Build the graded temporal mesh
\begin{equation*}
t_k=T\Bigl(\frac{k}{K}\Bigr)^{\gamma_g},\qquad
\tau_k=t_k-t_{k-1},\qquad
t_{k-\theta}=\theta t_{k-1}+(1-\theta)t_k.
\end{equation*}
\State Initialize $\tilde{v}^{0}(\boldsymbol{x})$ from the initial condition and set $\mathcal{V}^l(t_0;\boldsymbol{x})\gets 0$ for $1\le l\le N_{\exp}$.
\State Initialize the network parameters $\Theta^{0}$.

\For{$k=1$ to $K$}
  \State Warm start the trainable parameters by setting $\Theta\gets\Theta^{k-1}$.
  
  \For{$m=1$ to $m_{\mathrm{step}}$}
    \State Evaluate the current candidate
    \begin{equation*}
    \tilde{v}_{\Theta}^k(\boldsymbol{x})
    =\tilde{v}_{\Theta}(\boldsymbol{x},t_k),
    \qquad
    \tilde{v}_{\Theta}^{k-\theta}
    =\theta \tilde{v}^{k-1}+(1-\theta)\tilde{v}_{\Theta}^k.
     \end{equation*}
    
    \State Compute
    \begin{equation*}
    v_{\Theta}^{k-\theta}
    =\mathcal{N}\bigl[\tilde{v}_{\Theta}^{k-\theta};\lambda\bigr]
    -g\bigl(\tilde{v}_{\Theta}^{k-\theta}\bigr),
    \qquad
    \nabla_\tau \tilde{v}_{\Theta}^k
    =\tilde{v}_{\Theta}^k-\tilde{v}^{k-1}.
    \end{equation*}
    
    \State Update the candidate history variables $\mathcal{V}_\Theta^l(t_k;\cdot)$ using \eqref{recursive_1}, the frozen snapshots $\{\tilde{v}^j\}_{j=0}^{k-1}$, and the current increment $\nabla_\tau \tilde{v}_{\Theta}^k$.
    
    \State Compute the time-marching residual loss
    \begin{equation*}
    \mathcal{L}_M^k(\Theta)=\frac{1}{N_x}\sum_{i=1}^{N_x}
    \bigl|\mathcal{R}_\Theta^{k-\theta}(\boldsymbol{x}_f^i)\bigr|^2,
    \end{equation*}
    where $\mathcal{R}_\Theta^{k-\theta}$ is evaluated from \eqref{residual_1}.
    
    \State Update $\Theta$ by Adam using $\mathcal{L}_{M}^k(\Theta)$.
    
    \If{$\mathcal{L}_M^k(\Theta)<\varepsilon_{\mathrm{tol}}$}
       \State \textbf{break}
    \EndIf
  \EndFor
  
  \State Freeze the trained state:
  \begin{equation*}
  \Theta^{k}\gets\Theta,\qquad
  \tilde{v}^{k}\gets \tilde{v}_{\Theta^{k}}(\cdot,t_k),\qquad
  \mathcal{V}^l(t_k;\cdot)\gets \mathcal{V}_{\Theta^{k}}^l(t_k;\cdot),
  \quad 1\le l\le N_{\exp}.
  \end{equation*}
\EndFor

\State \Return $\{\tilde{v}^{k}\}_{k=1}^{K}$ and $\{\mathcal{V}^l(t_k;\cdot)\}_{k=1}^{K}$.
\end{algorithmic}
\end{algorithm}
%% this is the subsection related to time-marching scheme: HKD
\subsection{Time-Marching Configuration: Alikhanov-fPINN-M}
We now present the time-marching variant of the Alikhanov-fPINN framework, which we denote as Alikhanov-fPINN-M. 
This configuration sequentially enforces the SOE-accelerated second-order time-discrete residual along the graded temporal mesh. 
At the shifted time level $t_{k-\theta}$, the residual is defined as 
\begin{equation}\label{residual_1}
\begin{split}
\mathcal{R}^{k-\theta}_{\Theta}(\mathbf{x})
&:= \boldsymbol{\mathrm{a}}^{(0,k)}
\triangledown_\tau \tilde{v}^{k}_{\Theta}(\mathbf{x})
+ \sum_{l=1}^{N_{\exp}}
\nu^l \,\mathrm{e}^{-s^l\theta\tau_k}
\,\mathcal{V}^{l}_{\Theta}(t_k;\mathbf{x})  \\
&\qquad
+\mathcal{N}\bigl(\tilde{v}_{\Theta}(\mathbf{x},t_{k-\theta});\lambda\bigr)
-g(\mathbf{x},t_{k-\theta}),
\qquad 1\le k\le K .
\end{split}
\end{equation}
Here,
\begin{equation*}
  \triangledown_\tau \tilde{v}^{k}_{\Theta}(\mathbf{x})
:=\tilde{v}^{k}_{\Theta}(\mathbf{x})
-\tilde{v}^{k-1}(\mathbf{x}),
\end{equation*}
and the SOE history variables $\mathcal{V}^{l}_{\Theta}(t_k;\mathbf{x})$ are updated recursively according to \eqref{recursive_1}.
The Alikhanov-fPINN-M method imposes the discretization-consistent residual \eqref{residual_1} one time level at a time, so that the contribution of the temporal discretization can be examined transparently and reproducibly.
In this setting, previously computed solution snapshots are kept fixed and used as stored history when advancing to the next shifted time level.

Let
%\begin{equation*}
$0=t_0<t_1<\cdots<t_K=T$
%\end{equation*}
be the graded temporal mesh introduced in Section~\ref{graded_structure}. 
For each fixed $\mathbf{x}\in\Omega$, let $\tilde{v}^k(\mathbf{x})$ denote the frozen approximation at time $t_k$ produced by the marching procedure. 
For a candidate parameter vector $\Theta$, we write
\begin{equation}
\tilde{v}^{k}_{\Theta}(\mathbf{x})
:=\tilde{v}_{\Theta}(\mathbf{x},t_k),
\end{equation}
where $\tilde{v}_\Theta$ is the hard-constrained trial function~\eqref{ansatz}. 
The computation starts from the initial snapshot
\begin{equation}
\tilde{v}^{0}(\mathbf{x}) := 0,\qquad \mathbf{x}\in\Omega,
\end{equation}
which the ansatz \eqref{ansatz} satisfies exactly.
The history variables are initialized by
\begin{equation*}
   \mathcal{V}^{l}(t_0;\mathbf{x})=0,
   \qquad 1\le l\le N_q .
\end{equation*}

Assume that the snapshots $\{\tilde{v}^n(\mathbf{x})\}_{n=0}^{k-1}$ and the corresponding history states have already been computed and stored. 
At the next time level, the unknown snapshot is represented by the neural ansatz
\begin{equation*}
\tilde{v}_{\Theta}^{k}(\mathbf{x})
= \tilde{v}_\Theta(\mathbf{x},t_k)
= t_k\rho(\mathbf{x})v_{NN}(\mathbf{x},t_k;\Theta),
\qquad \mathbf{x}\in\Omega ,
\end{equation*}
while all earlier snapshots remain fixed. 
The shifted state used in the nonlinear and spatial terms is defined by
\begin{equation*}
\tilde{v}_{\Theta}^{k-\theta}(\mathbf{x})
:= \theta\,\tilde{v}^{k-1}(\mathbf{x})
+ (1-\theta)\tilde{v}_{\Theta}^{k}(\mathbf{x}),
\qquad \mathbf{x}\in\Omega .
\end{equation*}
The corresponding time increment is
\begin{equation*}
\triangledown_\tau \tilde{v}_{\Theta}^{k}(\mathbf{x})
:= \tilde{v}_{\Theta}^{k}(\mathbf{x})
- \tilde{v}^{k-1}(\mathbf{x}).
\end{equation*}

The nonlocal memory contribution in \eqref{residual_1} is evaluated through the SOE variables, which are updated recursively.
For $k\ge1$ and $1\le l\le N_q$, we set
\begin{equation*}
\mathcal{V}^{l}_{\Theta}(t_k;\mathbf{x})
= e^{-s^l\tau_k}\mathcal{V}^{l}(t_{k-1};\mathbf{x})
+ a^{(k,l)}\triangledown_\tau \tilde{v}_{\Theta}^{k}(\mathbf{x})
+ b^{(k,l)}
\bigl(\rho_k\triangledown_\tau \tilde{v}_{\Theta}^{k+1}(\mathbf{x})
-\triangledown_\tau \tilde{v}_{\Theta}^{k}(\mathbf{x})
\bigr),
\end{equation*}
with the initial value $\mathcal{V}^{l}(t_0;\mathbf{x})=0$. Using these quantities, the residual at the shifted level $t_{k-\theta}$ is evaluated through \eqref{residual_1}.

Let the set of interior spatial collocation points $\{\mathbf{x}_f^i\}_{i=1}^{N_x}\subset\Omega$ be the fixed points used at every marching step.
These points are generated once before training and remain unchanged throughout the sequential procedure. 
At step $k$, only the temporal evaluation level changes from one shifted time to the next. 
The stepwise loss functional is defined by
\begin{equation}\label{eq:loss_marching_step}
\mathcal{L}_{M}^{k}(\Theta)
:=\frac{1}{N_x}\sum_{i=1}^{N_x}
\bigl|\mathcal{R}_{\Theta}^{k-\theta}(\mathbf{x}_f^i)\bigr|^2,
\qquad 1\le k\le K ,
\end{equation}
and the trainable parameters at the $k$-th step are obtained from
\begin{equation}\label{eq:argmin_marching_step}
\Theta^{k}
:=\arg\min_{\Theta}\mathcal{L}_{M}^{k}(\Theta).
\end{equation}

After the minimization is completed, the newly computed snapshot is frozen as
\begin{equation*}
\tilde{v}^{k}(\mathbf{x})
:= \tilde{v}_{\Theta^{k}}^{k}(\mathbf{x}),
\qquad \mathbf{x}\in\Omega ,
\end{equation*}
and the memory variables are stored at the optimized parameters:
\begin{equation*}
\mathcal{V}^{l}(t_k;\mathbf{x})
:=\mathcal{V}_{\Theta^{k}}^{l}(t_k;\mathbf{x}),
\qquad 1\le l\le N_q,\quad \mathbf{x}\in\Omega .
\end{equation*}
These stored quantities are then used as the memory state for the subsequent marching step. In practice, the optimization at each new step is warm-started from the parameters obtained at the preceding step, which improves training stability on graded meshes and reduces fluctuations in the optimization error.

The Alikhanov-fPINN-M procedure sequentially advances from $k=1$ to $K$, producing the discrete-time surrogate $\{\tilde{v}^k\}_{k=1}^{K}$. 
Since the residual is imposed at each shifted level and the SOE memory variables are updated recursively, this configuration remains consistent with the accelerated second-order temporal discretization. 
In this work, the Alikhanov-fPINN-M method is used primarily as an auxiliary configuration to verify temporal convergence with respect to the maximum time step $\tau$.
By isolating the residual enforcement at each time level and treating the history terms with stored recursive variables, the method provides a controlled environment for evaluating the accuracy of the proposed time discretization.

%%%%%%%%%%%%%%%%%%%%%%%%%%%%%%%%%%%%%%%%%%%%%%%
\section{Numerical Results}\label{sec5}
This section presents the essential implementation details of the proposed numerical procedure. 
We apply our proposed Alikhanov-PINNs scheme to different test examples.
 To evaluate the efficiency and accuracy of the proposed numerical algorithm for various NfPDEs, we conduct forward and inverse problem tests, including cases involving the nonlinear subdiffusion equation, the generalized viscous Burgers equation, and cases without known exact solutions. 
Through analysis of network architecture, training data size, and adaptive variants of activation functions, we identify optimal Alikhanov-PINNs parameters. 
All computations are performed in Python on a workstation equipped with an Intel\textregistered\ Core\texttrademark\ i7-10870H CPU at 2.21 GHz and an NVIDIA GeForce GTX 1660 Ti GPU. 
Our implementation is based on TensorFlow~2.18, but can easily be adapted to other deep learning frameworks such as PyTorch.
 We compute the gradient with respect to the network parameters and inputs using the TensorFlow's built-in auto-differentiation module, and use the Adam optimizer to minimize the loss function. 
 
To quantify the accuracy and temporal convergence of the Alikhanov–PINN on nonuniform meshes, we use 
the maximum-norm error $E_{\infty}(K) = \max_k \|U^k - u^k\|_{\infty} $ 
and $L_2$-error $E_{2}(K)= \frac{ \|\tilde{v}(\cdot,t)-v(\cdot,t) \|_2}{ \|v(\cdot,t) \|_2}$. 
The temporal convergence rate % (TOC)
is then computed by the following formula:
\begin{equation}
    \operatorname{rate}_j\approx \log_2\Bigl(\frac{E_j(K)}{E_j(2K)}\Bigr),
\end{equation}
where $j=2$ or $j=\infty$.

%%%%%%%%%%%%%%%%%%%%%%%%%%%%%%%%%%%%%%%%% Section 3.1.
\subsection{Forward Problem:}

\begin{example}\label{ex1}

\textit{Nonlinear time-fractional subdiffusion equation (NTFSDE) with smooth exact solution}\\
In this example, we consider
\begin{align*}
{} _0^{\mathcal{C}}\partial_t^{\alpha} v+\mathcal{N}[v;\lambda]&=g(\mathbf{x},t),\\
v(\mathbf{x},0) &= 0,\quad \mathbf{x}\in \Omega,\\
v(\mathbf{x},t) &= 0,\quad \mathbf{x}\in \partial\Omega,
\end{align*}
%\begin{example}\label{ex1}
%\textit{Nonlinear time-fractional subdiffusion equation (NTFSDE) with smooth exact solution}\\
%In this example, we test the accuracy and consider 
%\begin{align*}
%{} _0^{\mathcal{C}}\mathcal{\partial}_t^{\alpha} v+\mathcal{N}[v;\lambda]&=0,\\
%v(\mathbf{x},0) &= 0,\quad \mathbf{x}\in \Omega,\\
%v(\mathbf{x},t) &= 0,\quad \mathbf{x}\in \partial\Omega,
%\end{align*}
%\Tmatthias{IC and BC Zero ? I.e. only the source $g$ prevents you from getting the trivial zero solution}
where $\mathcal{N}[v;\lambda]=\lambda_1 (v^3-v)-\lambda_2(v_{xx})+g(\cdot, t)$ with $(\lambda_1,\lambda_2)$ being unknown parameters which is the case of 1D problem setup. 
To extend this problem to 2D setup and 3D setup, we assume 
\begin{align*}
\mathcal{N}[v; \lambda]=\lambda_1(v^3-v)-\lambda_2(v_{xx}+v_{yy})-g(x,y,t),\\
\mathcal{N}[v; \lambda]=\lambda_1(v^3-v)-\lambda_2(v_{xx}+v_{yy}+v_{zz})-g(x,y,z,t)
\end{align*}
 with the smooth exact solutions being
 \begin{align*}
 v(\cdot,t) &= t^{2+\alpha} \sin(\pi x),\quad \Omega=[0,1],\\
 v(\cdot,t) &= t^{2+\alpha} \sin(\pi x)\sin(\pi y),\quad \Omega=[0,1]^2,\\
 v(\cdot,t) &= t^{2+\alpha} \sin(\pi x)\sin(\pi y)\sin(\pi z),\quad \Omega=[0,1]^3,
 \end{align*}
 and the corresponding body force $g$ is calculated based on these exact solutions.
 \end{example}
 The performance of the Alikhanov-PINNs is evaluated using the 1D NTFSDE problem to systematically examine the effects of activation functions (adaptive variant), neural network architecture and the size of the training data sample. 
 The solution is represented as $\tilde{v}=t \mathbf{x}(1-\mathbf{x}) v_{NN}(\mathbf{x},t,\Theta)$, $\tilde{v}$ is the Alikhanov-PINNs solution and IBC are imposed via hard constraints. 
 The training data are sampled from a nonuniform mesh ($\gamma=1$). 
 
\begin{figure}[htbp]
  \centering
  \begin{subfigure}[b]{0.48\textwidth}
    \includegraphics[width=\textwidth]{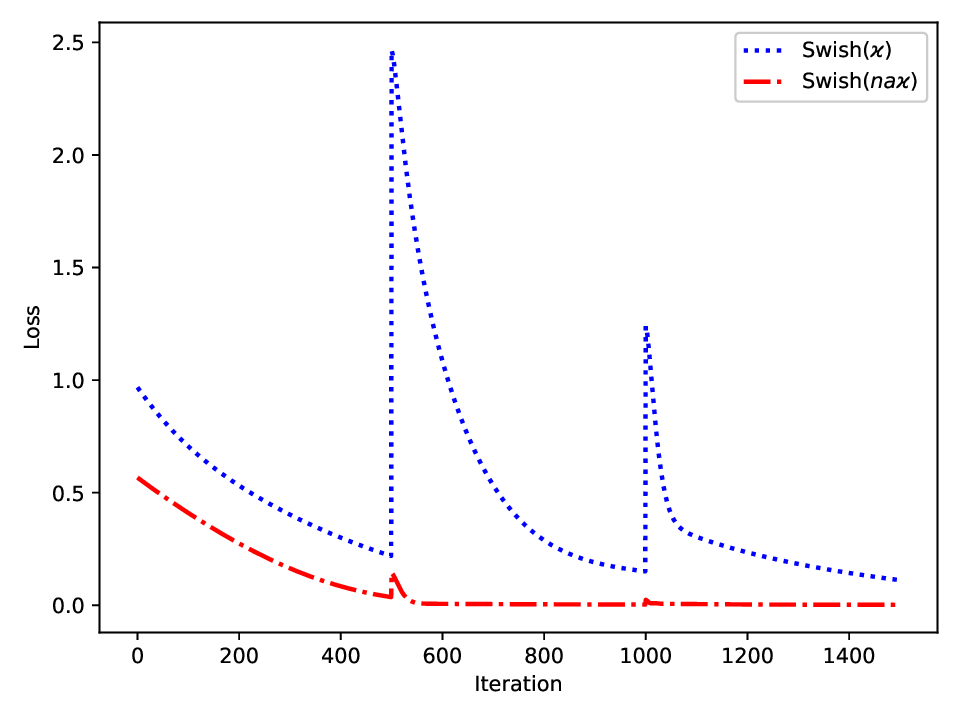}
    \caption{}
    \label{fig:loss_swish_vs_adap}
  \end{subfigure}
  \hfill
  \begin{subfigure}[b]{0.48\textwidth}
    \includegraphics[width=\textwidth]{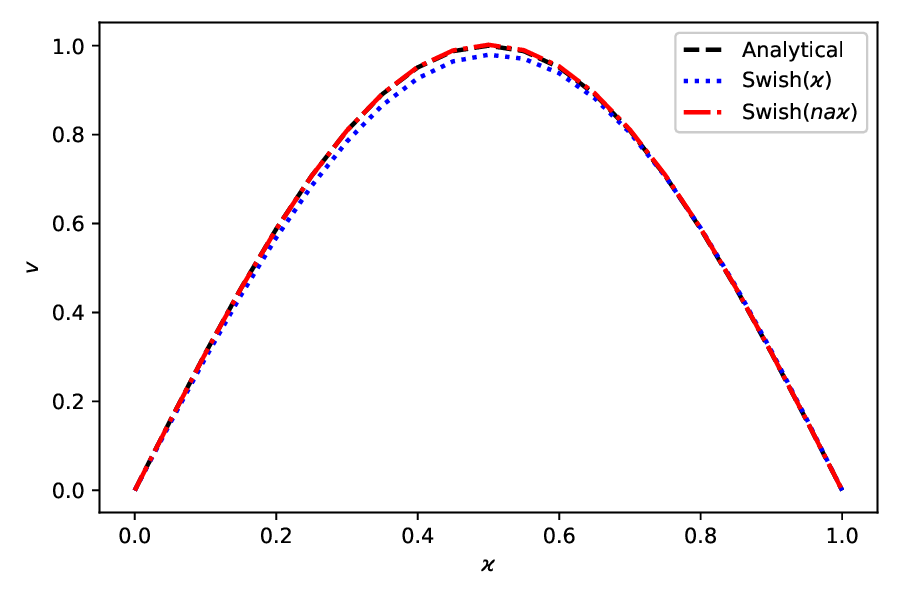}
    \caption{}
    \label{fig:swish_vs_adap}
  \end{subfigure}
  \caption{Comparison of loss evolution and Alikhanov--XfPINNs approximation results for $\Swish(x)$ and $\Swish(na x)$ ($\alpha=0.5$): (a) training loss progression; (b) exact and predicted field profiles.}
  \label{fig:swish_adap_comparison}
\end{figure}

 To determine the most effective activation function, we compare six adaptive variants in Table~\ref{tab_six_adap_f} for $\alpha=0.5$. 
The $\Swish(na\mathbf{x})$ function yields the lowest $\mathcal{E}_{\infty}$-error and relative $\mathcal{E}_2$-error outperforming its fixed counterpart $\Swish(\mathbf{x})$ by capturing the nonlinear dynamics more accurately. Figure~\ref{fig:swish_adap_comparison} confirms that $\Swish(na\mathbf{x})$ leads to a faster loss decay and better agreement with the analytical solution, highlighting the importance of adaptive activations in improving the training efficiency. 
 
\begin{figure}[htbp]
  \centering
  \begin{subfigure}[b]{0.48\textwidth}
    \includegraphics[width=\textwidth]{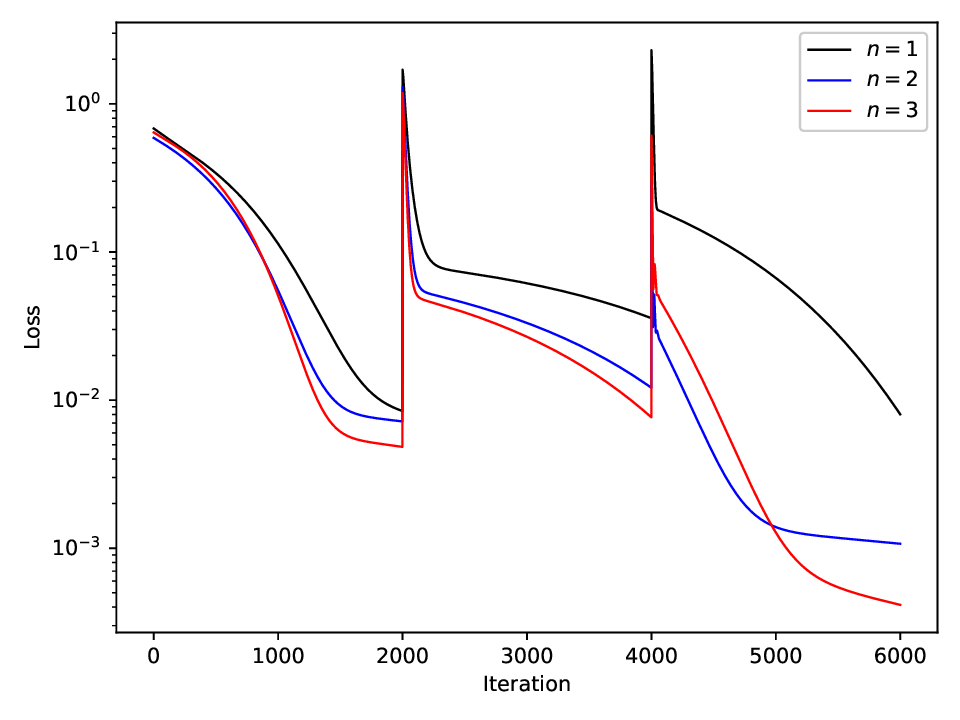}
    \caption{}
    \label{fig:diff_loss}
  \end{subfigure}
  \hfill
  \begin{subfigure}[b]{0.48\textwidth}
    \includegraphics[width=\textwidth]{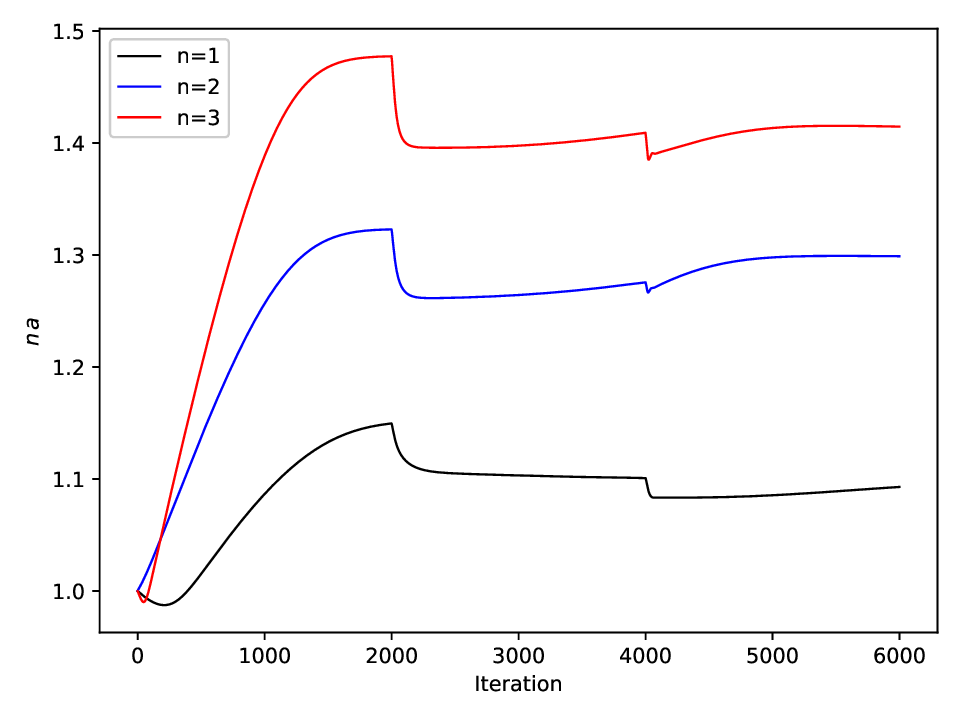}
    \caption{}
    \label{fig:diff_loss_na}
  \end{subfigure}
  \caption{Evolution of the training loss and $na$ across iterations with $\alpha=0.5$ for adaptive activation functions with various $n$ : (a) loss trajectory; (b) $na$.}
  \label{fig:loss_diff_na_comparison}
\end{figure}

\begin{figure}[htbp]
  \centering
  \begin{subfigure}[b]{0.32\textwidth}
    \includegraphics[width=\textwidth]{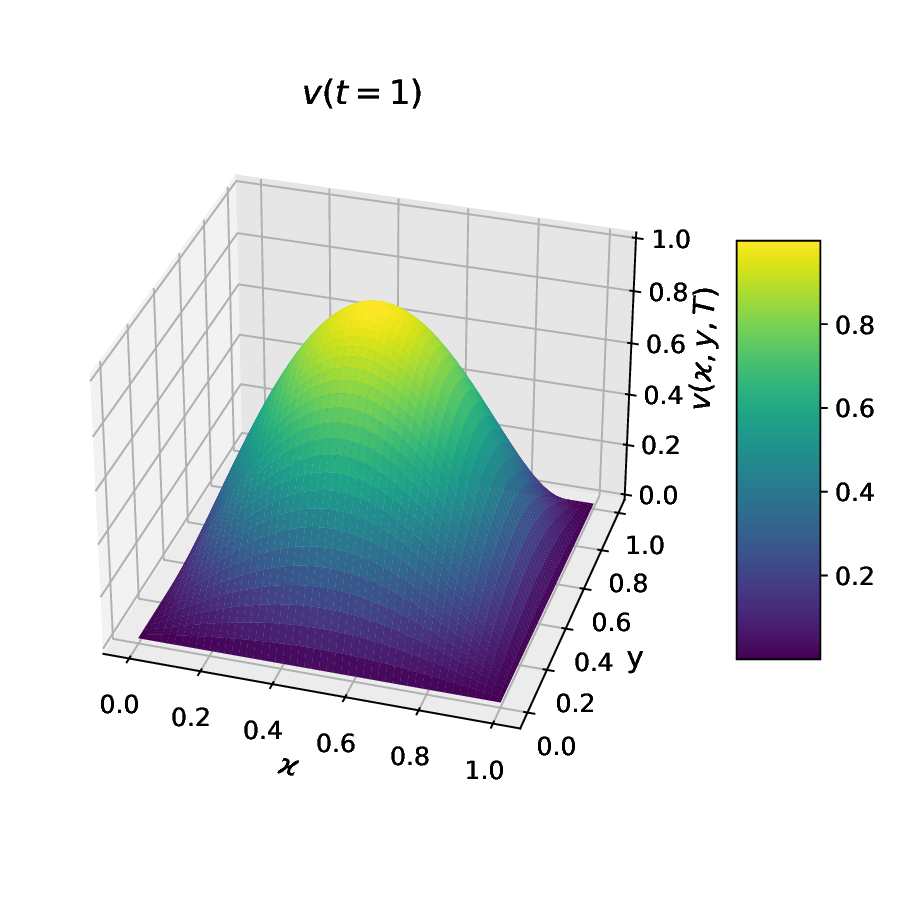}
    \caption{}
    \label{fig:ex1_exact}
  \end{subfigure}
  \hfill
  \begin{subfigure}[b]{0.32\textwidth}
    \includegraphics[width=\textwidth]{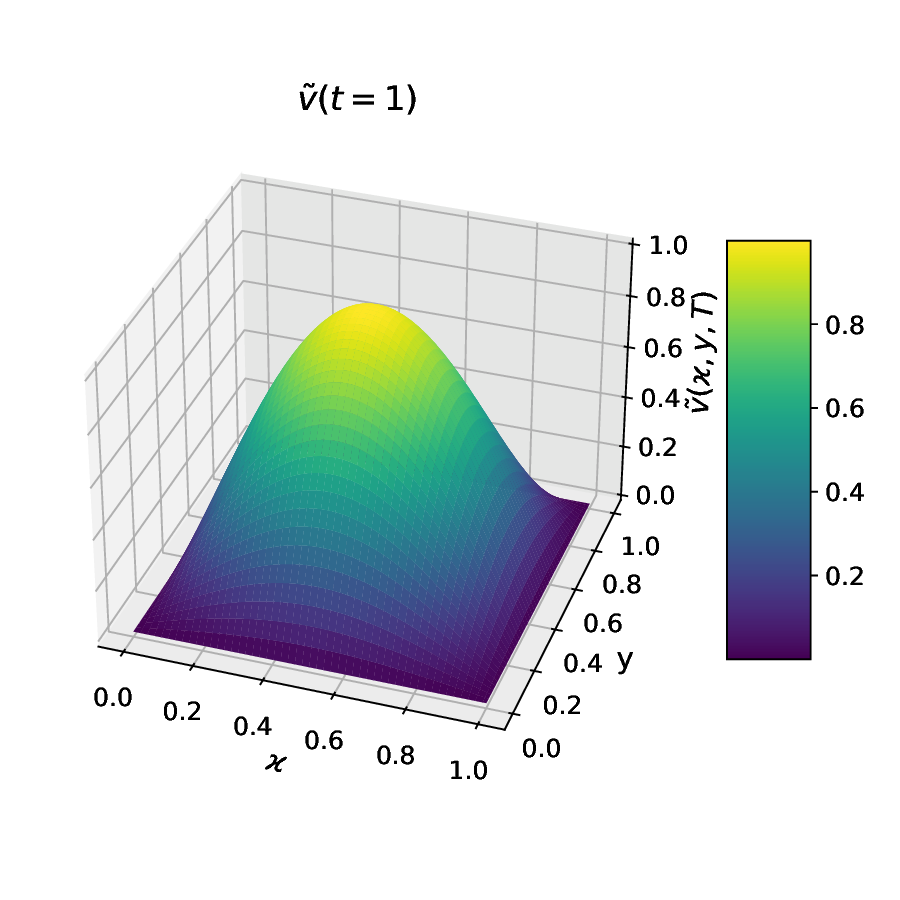}
    \caption{}
    \label{fig:ex1_num}
  \end{subfigure}
  \hfill
  \begin{subfigure}[b]{0.32\textwidth}
    \includegraphics[width=\textwidth]{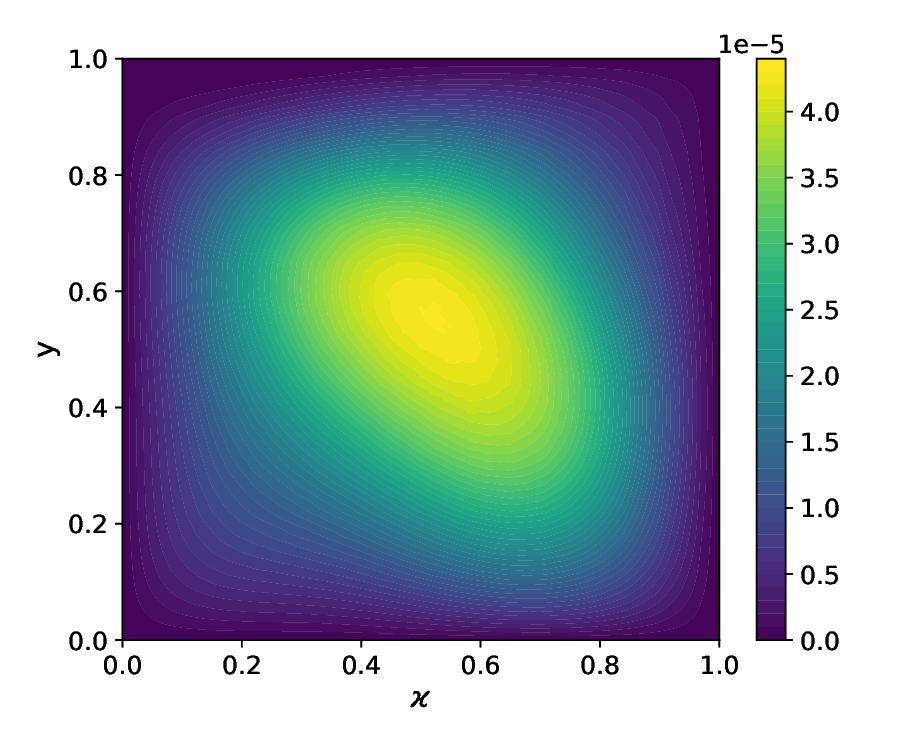}
    \caption{}
    \label{fig:ex1_err}
  \end{subfigure}
  \caption{Comparison of 2D-NTFSDe with smooth exact solution at $\alpha=0.9$ and $K=8$: (a) exact solution $v$, (b) Alikhanov--XfPINNs approximation $\tilde{v}$, and (c) absolute error $|v-\tilde{v}|$.}
  \label{fig:ex1_comparison}
\end{figure}

 Figures~\ref{fig:diff_loss} and \ref{fig:diff_loss_na} analyze the influence of the scale factor $n$ on training dynamics. 
 A higher $n$ value accelerates the convergence due to more aggressive parameter updates; however, overly large values may cause instability. 
 Convergence of the product $na$ during training indicates stabilization of the adaptive training rate. 
 %The effect of the fully connected neural network depth and width on accuracy is  shown in figure \ref{fig9a}-\ref{fig9b}. Deeper and wider networks enhance the representational capacity and reduce the prediction error. However, excessive depth or width can lead to overfitting and increased training cost due to the rise in trainable parameters.
 %
 %Figures \ref{fig10a}-\ref{fig10b} reveal the variation of relative $\mathcal{E}_2$-error of Alikhanov-PINN with number of spatial and temporal training points for $\alpha=0.3,0.6,0.9$, which demonstrate that increasing temporal training points significantly improves the solution accuracy while spatial refinement yields diminishing returns. Beyond a saturation point, further temporal refinement offers negligible accuracy gain. XfPINNs, employing domain decomposition, divide the computational domain into subdomains, each modeled by a separate neural network.This reduces model complexity per subdomain, improves generalization and convergence, and enables parallel training for enhanced computational efficiency 
 %
 Based on these results, $\Swish(na x)(n = 3)$ is selected as the optimal activation function.”. 
 Table~\ref{tab:layer_1} summarizes the final network configuration and associated training times for the 1D-3D problems. {Table~\ref{tab:example1} reports the $\mathcal{E}_{\infty}(K)$-errors and the corresponding temporal convergence rates obtained by the Alikhanov-fPINN-M for different values of $\alpha$. The results show that, when the grading parameter is chosen as $\gamma=2/\alpha$, the Alikhanov-fPINN-M recovers the expected second-order accuracy in the temporal direction. Table~\ref{tab:ex1_adaptive_training} reports the adaptive-training accuracy for the NTFSDE with smooth exact solution in Example~\ref{ex1}. The results show that Alikhanov-XfPINNs consistently produce smaller $\mathcal{E}_{\infty}(K)$ and $\mathcal{E}_{2}(K)$ errors than Alikhanov-PINNs for all tested fractional orders, which confirms that the adaptive XfPINN formulation improves the approximation accuracy under the same training setting.} Figure~\ref{fig:ex1_comparison} displays the exact solution, the corresponding numerical approximation, and the absolute error obtained from the proposed scheme.
 %Figure~\ref{fig:ex1_comparison} displays the exact solution, the corresponding numerical approximation, and the absolute error obtained from the proposed scheme.

 %Figure  \ref{fig:ex1_comparison} illustrates the close agreement between the predicted and analytical solution with errors uniformly small even in higher dimensions. Table \ref{tab6} contrasts the computational complexity of the traditional FDM and Alikhanov-PINNs. FDM scales as $\mathcal{O}(M_\varkappa M_y M_\varkappa N_1^2)$, rapidly growing with dimensionality. In contrast, Alikhanov-PINNs achieves comparable accuracy with fewer training points and reduced complexity $\mathcal{O}(J N_wN_2\mathcal{N}_q)$ demonstrating its superior efficiency in high-dimensional scenarios.
 
\begin{table}[htbp]
\scriptsize
\centering
\caption{Maximum‐norm error $\mathcal{E}_{\infty}(K)$ and relative $L^2$-error $\mathcal{E}_{2}(K)$ of the Alikhanov--XfPINNs with adaptive variant ($\alpha=0.5$, $K=16$).}
\begin{tabular}{lcc}
\hline
Adaptive Activation Function & $\mathcal{E}_{2}(K)$ & $\mathcal{E}_{\infty}(K)$ \\
\hline
$\operatorname{Sigmoid}(\mathbf{x})$           & $3.975e-04$ & $3.298e-04$ \\
$\Swish(\mathbf{x})$             & $2.685e-04$ & $1.678e-04$ \\
$\operatorname{ReLU}(\mathbf{x})$              & $9.479e-04$ & $9.893e-04$ \\
$\tanh(\mathbf{x})$              & $3.478e-04$ & $2.523e-04$ \\
$\operatorname{SeLU}(\mathbf{x})$              & $8.577e-04$ & $1.178e-03$ \\
$\mathbf{x}\,\tanh(\mathbf{x})$   & $3.012e-04$ & $2.223e-04$ \\
\hline
\end{tabular}
\label{tab_six_adap_f}
\end{table}

%%%%%%%%%%%%%%%%%%%%%%%%%%%%%%%%%%%%
\begin{table}[htbp]
\scriptsize
\caption{Overview of Alikhanov--XfPINNs architecture details and the quantity of training samples utilized.}
\centering
\renewcommand{\arraystretch}{1.2} % increase row height if you like
\setlength{\tabcolsep}{0pt}       % we'll distribute spacing manually
\begin{tabular*}{\textwidth}{@{\extracolsep{\fill}} c c c c c }
\hline
$K$ & $K_{\mathrm{total}}$ & $L$ & $N_i$ & $P$ \\
\hline
8   & 288  & 2 & 14 & 281 \\
16  & 576  & 3 & 15 & 556 \\
32  & 1152  & 4 & 18 & 1117 \\
\hline
\end{tabular*}
\label{tab:layer_1}
\end{table}

%%%%%%%%%%%%%%%%%%%%%%%%

%%%%%%%%%%%%%%%%%%%%%%%%%%%%%%%%%%%%%%%%%%
\begin{table}[htbp]
\centering
\scriptsize
\caption{$\mathcal{E}_{\infty}(K)$-errors and temporal order of convergence (TOC) rates of the Alikhanov--fPINN-M time-marching formulation for the 2D-NTFSDE with smooth exact solution in Example~\ref{ex1}.} 
\begin{tabular}{c c c c}    
\hline
$\alpha$ 
& $K$ 
& \multicolumn{2}{c}{Alikhanov-fPINN-M} \\
\hline
& & $\mathcal{E}_{\infty}(K)$ 
& $rate_{\tau}$ \\
\hline
$0.30$ & $2^3$ & $2.582e-04$ & $--$ \\
       & $2^4$ & $6.730e-05$ & $1.939$ \\
       & $2^5$ & $1.715e-05$ & $1.972$ \\
\textit{TOC} &  &  & $2.0$ \\
\hline
$0.60$ & $2^3$ & $4.651e-04$ & $--$ \\
       & $2^4$ & $1.228e-04$ & $1.921$ \\
       & $2^5$ & $3.147e-05$ & $1.964$ \\
\textit{TOC} &  &  & $2.0$ \\
\hline
$0.90$ & $2^3$ & $6.521e-04$ & $--$ \\
       & $2^4$ & $1.730e-04$ & $1.914$ \\
       & $2^5$ & $4.399e-05$ & $1.975$ \\
\textit{TOC} &  &  & $2.0$ \\
\hline
\end{tabular}
\label{tab:example1}
\end{table}

\begin{table}[htbp]
\centering
\scriptsize
%\captionsetup{labelfont={color=red}, textfont={color=red}}
%\color{red}
\caption{Maximum-norm error $\mathcal{E}_{\infty}(K)$ and relative $L^2$ error $\mathcal{E}_{2}(K)$ obtained by the adaptive-training Alikhanov-PINNs and Alikhanov-XfPINNs for the NTFSDE with smooth exact solution in Example~\ref{ex1}.}
\begin{tabular}{c c c c c}
\hline
$\alpha$ 
& \multicolumn{2}{c}{Alikhanov-PINNs} 
& \multicolumn{2}{c}{Alikhanov-XfPINNs} \\
\cline{2-5}
& $\mathcal{E}_{\infty}(K)$ 
& $\mathcal{E}_{2}(K)$ 
& $\mathcal{E}_{\infty}(K)$ 
& $\mathcal{E}_{2}(K)$ \\
\hline
$0.30$ & $3.435e-03$ & $4.536e-03$ & $1.396e-03$ & $3.589e-03$ \\
$0.60$ & $2.185e-03$ & $3.874e-03$ & $1.057e-03$ & $1.247e-03$ \\
$0.90$ & $1.864e-03$ & $3.318e-03$ & $9.735e-04$ & $1.196e-03$ \\
\hline
\end{tabular}
\label{tab:ex1_adaptive_training}
\end{table}

\begin{figure}
  \centering
  \makebox[\textwidth][c]{%
  \includegraphics[width=12.9cm, height=6.75cm]{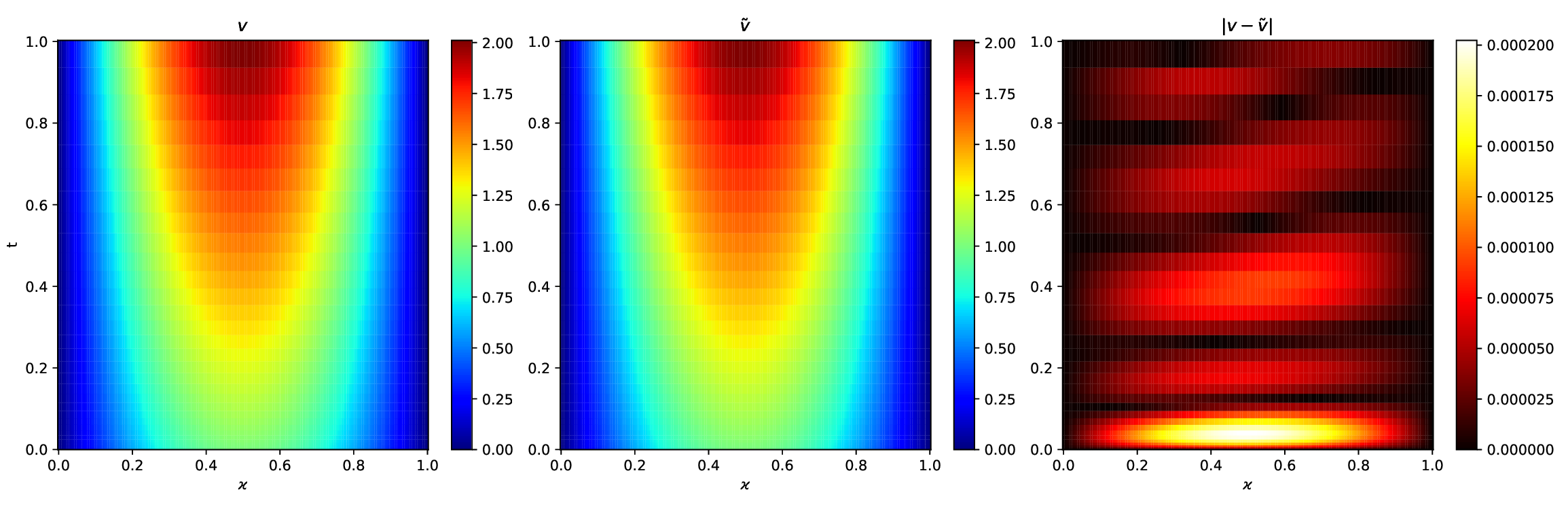}%
  }
  \caption{Comparison of exact solution $v$, Alikhanov--XfPINNs prediction $\tilde{v}$ and absolute error $|v-\tilde{v}|$ for time fractional generalized Burgers' problem.}
  \label{example2}
\end{figure}

%%%%%%%%%%%%%%%%%%%%%%%%%%%%%%%
\begin{example}\textit{Generalized Burgers' equation with nonsmooth exact solution}\\
To verify the effectiveness of our nonuniform scheme, we consider the time-fractional generalized viscous Burgers' problem, for which the exact solution exhibits an initial singularity at $t=0$. 
We consider 
\begin{equation}
{} _{0}^{\mathcal{C}}\mathcal{\partial}_t^{\alpha}v+\mathcal{N}[v;\lambda] = g(x,t),
\end{equation}
where $\mathcal{N}[v;\lambda] = \lambda_1v^p v_x+\lambda_2 v_{xx}$, $p>0$ and the exact solution is $v= (1+\omega_{1+\alpha}(t))\sin(\pi x)\sin(\pi y)$, $(\mathbf{x},t)\in[0,1]\times[0,1]$, for the initial condition $v_0(x)=\sin(\pi x)$. 
\end{example}
Due to the lack of first order time continuity, the exact solution exhibit singular behavior as $t\to 0$. 
To address this we utilize the same fully connected neural network(FCNN) structure as in Example~\ref{ex1}, embedding the initial and boundary conditions as hard constraint. 
The network approximation is expressed as $\tilde{v}(x,t)=tx(1-x)v_{\text{NN}}(x,t,\Theta)$. 
To resolve the initial singularity, we examine the effect of graded temporal meshes with refinement parameters $\gamma=1$(uniform), $\gamma=\frac{2}{\alpha},\frac{3-\alpha}{\alpha}$.
Figure~\ref{example2} presents the absolute error profiles when $\alpha=0.9$. 
Tables~\ref{tab:example21}, \ref{tab:example22} and \ref{tab:example23} present the $\mathcal{E}_{\infty}$-error obtained from the Alikhanov-PINNs and Alikhanov--XfPINNs trained on a different grading parameter.
As can be seen from the tables, the optimal grading parameter that offers the best convergence rate is clearly $\gamma=\frac{2}{\alpha}$. {Table~\ref{tab:ex2_adaptive_training} reports the adaptive-training accuracy of Alikhanov-PINNs and Alikhanov--XfPINNs for the generalized viscous Burgers equation with nonsmooth exact solution. The results show that the adaptive Alikhanov-XfPINNs consistently reduce both $\mathcal{E}_{\infty}(K)$ and $\mathcal{E}_{2}(K)$ for all tested fractional orders, confirming the effectiveness of the adaptive XfPINN formulation for this weakly singular benchmark.}

%%%%%%%%%%%%%%%%%%%%%%%%%%%%
\begin{table}[htbp]
\centering
\scriptsize
\caption{$\mathcal{E}_{\infty}$–errors and temporal convergence rates of Alikhanov-fPINN-M for the generalized viscous Burgers equation with nonsmooth exact solution ($\alpha = 0.3$).}
\begin{tabular}{c c c c c c c}    
\hline
   $\alpha$ & $\gamma$ & $K$ & \multicolumn{2}{c}{$p=2$} & \multicolumn{2}{c}{$p=4$} \\
\hline
 & & & $\mathcal{E}_{\infty}(K)$ & $rate_{\tau}$ & $\mathcal{E}_{\infty}(K)$ & $rate_{\tau}$ \\
\hline
&&&\\
$0.30$ & $1$             & $2^3$ & $8.122e-03$  & --      & $8.095e-03$  & --     \\
       &                 & $2^4$ & $7.472e-03$  & $0.120$ & $6.103e-03$  & $0.409$\\
       &                 & $2^5$ & $6.670e-03$  & $0.164$ & $4.350e-03$  & $0.489$\\
\it{TOC}&                &       &             & $-$   &              & $-$  \\
\hline
&&&\\
       & $2/\alpha$      & $2^3$ & $2.812e-03$  & --      & $2.700e-03$  & --     \\
       &                 & $2^4$ & $8.315e-04$  & $1.758$ & $7.514e-04$  & $1.845$\\
       &                 & $2^5$ & $2.335e-04$  & $1.832$ & $2.030e-04$  & $1.888$\\
\it{TOC}&                &       &             & $2.0$   &              & $2.0$  \\
\hline
&&&\\
       & $(3-\alpha)/\alpha$ & $2^3$ & $1.845e-03$ & --      & $2.861e-03$  & --     \\
       &                 & $2^4$ & $5.094e-04$  & $1.857$ & $8.544e-04$  & $1.744$\\
       &                 & $2^5$ & $1.362e-04$  & $1.903$ & $2.315e-04$  & $1.884$\\
\it{TOC}&                &       &             & $2.0$   &              & $2.0$  \\
\hline
\end{tabular}
\label{tab:example21}
\end{table}

%%%%%%%%%%%%%%%%%%%%%%%%%%%%
\begin{table}[htbp]
\centering
\scriptsize
\caption{$\mathcal{E}_{\infty}$–errors and temporal convergence rates of Alikhanov-fPINN-M for the generalized viscous Burgers equation with nonsmooth exact solution ($\alpha = 0.6$).}
\begin{tabular}{c c c c c c c}    
\hline
   $\alpha$ & $\gamma$ & $K$ & \multicolumn{2}{c}{$p=2$} & \multicolumn{2}{c}{$p=4$} \\
\hline
 & & & $\mathcal{E}_{\infty}(K)$ & $rate_{\tau}$ & $\mathcal{E}_{\infty}(K)$ & $rate_{\tau}$ \\
\hline
&&&\\
$0.60$ & $1$             & $2^3$ & $4.439e-03$  & --      & $4.558e-03$  & --     \\
       &                 & $2^4$ & $3.779e-03$  & $0.232$ & $3.635e-03$  & $0.326$\\
       &                 & $2^5$ & $2.825e-03$  & $0.419$ & $2.593e-03$  & $0.487$\\
\it{TOC} &                &       &             & $-$     &              & $-$    \\
\hline
&&&\\
       & $2/\alpha$      & $2^3$ & $1.295e-03$  & --      & $1.437e-03$  & --     \\
       &                 & $2^4$ & $3.811e-04$  & $1.765$ & $5.037e-04$  & $1.512$\\
       &                 & $2^5$ & $1.042e-04$  & $1.870$ & $1.357e-04$  & $1.892$\\
\it{TOC} &                &       &             & $2.0$   &              & $2.0$  \\
\hline
&&&\\
       & $(3-\alpha)/\alpha$ & $2^3$ & $1.047e-03$ & --      & $1.226e-03$  & --     \\
       &                 & $2^4$ & $2.970e-04$  & $1.818$ & $3.463e-04$  & $1.824$\\
       &                 & $2^5$ & $7.801e-05$  & $1.929$ & $9.405e-05$  & $1.881$\\
\it{TOC} &                &       &             & $2.0$   &              & $2.0$  \\
\hline
\end{tabular}
\label{tab:example22}
\end{table}

%%%%%%%%%%%%%%%%%%%%%%%%%%%%%%%%%%%%%%%%%%
\begin{table}[htbp]
\centering
\scriptsize
\caption{$\mathcal{E}_{\infty}$–errors and temporal convergence rates of Alikhanov-fPINN-M for the generalized Burgers equation with nonsmooth exact solution ($\alpha = 0.9$).}
\begin{tabular}{c c c c c c c}    
\hline
   $\alpha$ & $\gamma$ & $K$ & \multicolumn{2}{c}{$p=2$} & \multicolumn{2}{c}{$p=4$} \\
\hline
 & & & $\mathcal{E}_{\infty}(K)$ & $rate_{\tau}$ & $\mathcal{E}_{\infty}(K)$ & $rate_{\tau}$ \\
\hline
&&&\\
$0.90$ & $1$             & $2^3$ & $2.069e-03$  & --      & $1.895e-03$  & --     \\
       &                 & $2^4$ & $1.573e-03$  & $0.395$ & $1.475e-03$  & $0.361$\\
       &                 & $2^5$ & $1.116e-03$  & $0.495$ & $9.390e-04$  & $0.652$\\
\it{TOC}&                &       &             & $-$     &              & $-$    \\
\hline
&&&\\
       & $2/\alpha$      & $2^3$ & $5.773e-04$  & --      & $6.631e-04$  & --     \\
       &                 & $2^4$ & $1.516e-04$  & $1.929$ & $2.302e-04$  & $1.526$\\
       &                 & $2^5$ & $3.898e-05$  & $1.959$ & $6.019e-05$  & $1.935$\\
\it{TOC}&                &       &             & $2.0$   &              & $2.0$  \\
\hline
&&&\\
       & $(3-\alpha)/\alpha$ & $2^3$ & $5.210e-04$ & --      & $5.934e-04$  & --     \\
       &                 & $2^4$ & $1.345e-04$  & $1.954$ & $1.822e-04$  & $1.703$\\
       &                 & $2^5$ & $3.418e-05$  & $1.976$ & $4.863e-05$  & $1.906$\\
\it{TOC}&                &       &             & $2.0$   &              & $2.0$  \\
\hline
\end{tabular}
\label{tab:example23}
\end{table}

%\end{example}
%%%%%%%%%%%%%%%%%%%%%%%%%%%%%%%%%%%% Table 9
\begin{table}[htbp]
\centering
\scriptsize
%\captionsetup{labelfont={color=red}, textfont={color=red}}
%\color{red}
\caption{Maximum-norm error $\mathcal{E}_{\infty}(K)$ and relative $L^2$ error $\mathcal{E}_{2}(K)$ obtained by the adaptive-training Alikhanov-PINNs and Alikhanov-XfPINNs for the generalized viscous Burgers equation with nonsmooth exact solution.}
\begin{tabular}{c c c c c}
\hline
$\alpha$
& \multicolumn{2}{c}{Alikhanov-PINNs}
& \multicolumn{2}{c}{Alikhanov-XfPINNs} \\
\cline{2-5}
& $\mathcal{E}_{\infty}(K)$
& $\mathcal{E}_{2}(K)$
& $\mathcal{E}_{\infty}(K)$
& $\mathcal{E}_{2}(K)$ \\
\hline
$0.30$ & $8.024e-03$ & $8.474e-03$ & $7.575e-03$ & $7.834e-03$ \\
$0.60$ & $6.825e-03$ & $7.284e-03$ & $5.927e-03$ & $6.135e-03$ \\
$0.90$ & $6.015e-03$ & $6.275e-03$ & $5.285e-04$ & $5.593e-03$ \\
\hline
\end{tabular}
\label{tab:ex2_adaptive_training}
\end{table}

%\begin{example}
%KdV equation

%\end{example}

%%%%%%%%%%%%%%%%%%%%%%%%%%%%%%% EXAMPLE
\begin{example}
To examine a more realistic setting, we consider a problem without a closed-form exact solution, prescribed only by the initial condition and with no source term
\begin{equation}\label{unknown_prb}
{} _{0}^{\mathcal{C}}\mathcal{\partial}_t^{\alpha} v(\cdot,t)+\mathcal{N}[v;\lambda]=0,
\end{equation}
where $\mathcal{N}[v;\lambda]=\lambda \cdot \nabla v+v(v-1)(\rho-v)$ with $\lambda_1=\lambda_2=1, \rho>1$ and $\mathbf{x}\in\Omega, t\in [0,T]$. 

The exact solution to \eqref{unknown_prb} remains unknown. 
However, for the special case $\alpha = 1$, Table~\ref{tab:analytical_solutions} provides analytical expressions corresponding to the classical integer-order Fisher–Nagumo (FN) equation in one and two dimensions. 
These expressions are used to approximate the initial and boundary conditions (IBCs).

To solve the time-fractional Fisher–Nagumo (TFFN) problem with $\rho = 1$, and $\phi = \pi/2$, we use the Alikhanov--XfPINNs framework incorporating soft-constrained IBCs.
The predicted solution is denoted by $\tilde{v}(\mathbf{x},t) = v_{\text{NN}}(\mathbf{x},t;\Theta)$, with training data sampled from random boundary points, initial conditions, and a uniform spatial grid. 
Table~\ref{tab:tffn_architecture} summarizes the optimal network architecture and its associated training epoch.

Given the absence of an analytical solution for the TFFN model, we train the Alikhanov--XfPINN at fractional orders $\alpha\in\{0.30,0.60,0.90\}$ and assess its ability to recover the integer‑order FN solution $(\alpha=1)$. 
As illustrated in Figure~\ref{fig:TFFN}, the neural approximation increasingly aligns with the known solution as $\alpha \to 1$, validating the effectiveness of Alikhanov--XfPINN for inverse problems where the exact solution is unknown.

%%%%%%%%%%%%%%%%%%%%%%%%%%%%%%%%%%%
\begin{table}[htbp]
\centering
\scriptsize
\caption{Integer‑order FN $(\alpha=1)$: exact solutions $v(\mathbf{x},t)$ in $1D$ and $2D$.}

\begin{tabular}{l l}
\hline
Domain & $v(x,t)$ \\
\hline
$\Omega = [0,1]$ 
  & $\displaystyle 1/2 + 1/2\,\tanh\Bigl(\frac{1}{\sqrt{8}}[x-t\frac{2\rho-1}{\sqrt{2}} ]  \Bigr)$ \\[1ex]
$\Omega = [0,1]^2$ 
  & $\displaystyle 1/2 + 1/2\,\tanh\Bigl(\frac{1}{\sqrt{8}}[x \sin \varphi +y\cos\varphi]- t\frac{2\rho-1}{\sqrt{2}} \Bigr)$ \\
\hline
\end{tabular}
\label{tab:analytical_solutions}
\end{table}

%%%%%%%%%%%%%%%%%%%%%%%%%%%%%%%%%%%%%%%%
\begin{table}[htbp]
\scriptsize
\centering
\caption{Network architectures and training parameters for 1D‐TFFN and 2D‐TFFN models.}
\begin{tabular}{l c c l c c}
\hline
Equation & Depth & Width & Optimizer & Learning rate & Epoch \\
\hline
1D‐TFFN & $1$ & $30$ & Adam & $1e-04$ & $1500$ \\
2D‐TFFN & $3$ & $14$ & Adam & $1e-04$ & $1000$ \\
\hline
\end{tabular}
\label{tab:tffn_architecture}
\end{table}

%%%%%%%%%%%%%%%%%%%%%%%%%%%%%%%%%
\begin{figure}[htbp]
  \centering
  \begin{subfigure}[b]{0.48\textwidth}
    \includegraphics[width=\textwidth]{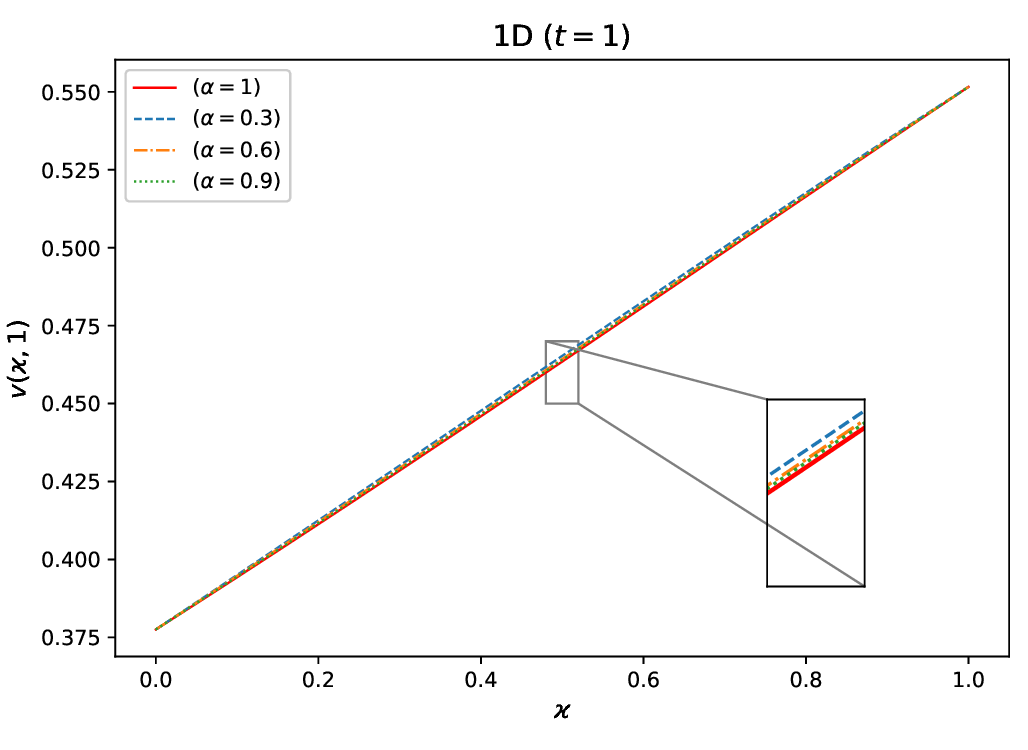}
    \caption{}
    \label{fig:TFFN_1D}
  \end{subfigure}
  \hfill
  \begin{subfigure}[b]{0.48\textwidth}
    \includegraphics[width=\textwidth]{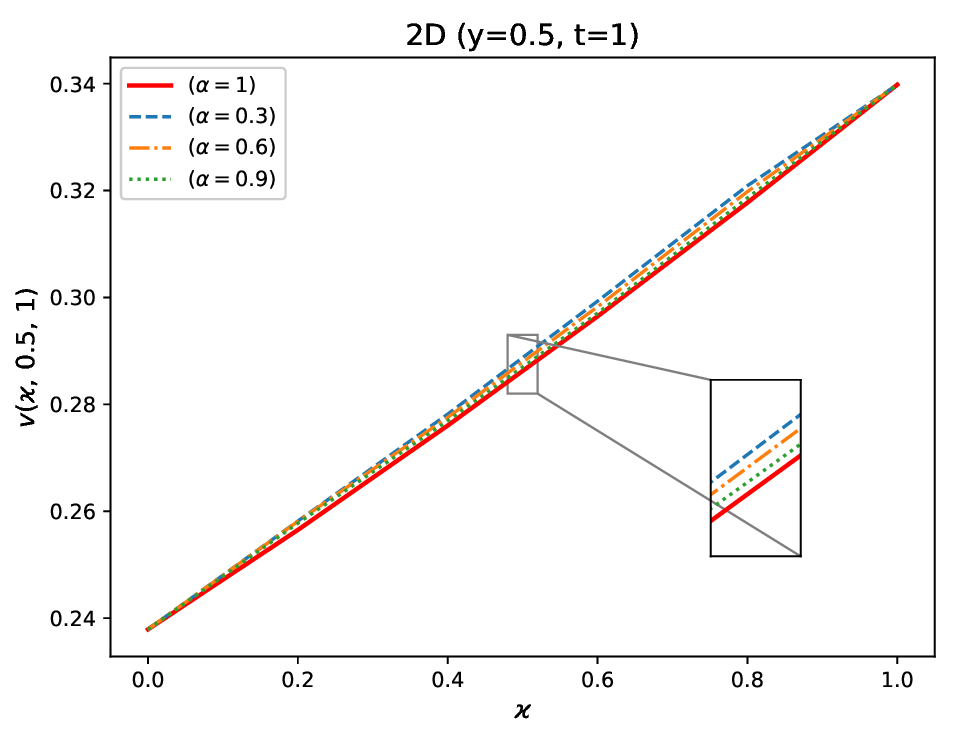}
    \caption{}
    \label{fig:TFFN_2D}
  \end{subfigure}
  \caption{Alikhanov--XfPINN predictions for $\alpha=0.3, 0.6, 0.9$ vs.\ exact integer-order FN solutions: (a) $1D$-TFFN; (b) $2D$-TFFN.}
  \label{fig:TFFN}
\end{figure}
\end{example}

%%%%%%%%%%%%%%%%%%%%%%%%%%%%%%%%%%%%%%%%%%%%%%%% Section 3.2.
\subsection{Inverse Problem}
\begin{example}
Here, we consider the inverse problem of NfPDEs with constant coefficients,
\begin{equation}
{} _{0}^{\mathcal{C}}\mathcal{\partial}_t^{\alpha}v+\mathcal{N}[v;\lambda]=g(\mathbf{x},t),
\end{equation}
with $\mathcal{N}[v;\lambda]=\lambda_1(v_{xx}+v_{yy})+\lambda_2 f(v) $ and
\begin{align*}
v(\mathbf{x},0)&=v_1(\mathbf{x}), \,\mathbf{x}\in\Omega,\quad
v(\mathbf{x},t)=v_2(\mathbf{x}), \,\mathbf{x,t}\in \partial\Omega, \,t\in(0,T),\quad 
v(\mathbf{x},T)=v_3(\mathbf{x},t), \,\mathbf{x}\in \Omega,
\end{align*}
where $\alpha,\lambda_1,\lambda_2$ being the unknown parameters. Given the reaction term $f(v)$, 
source term $g(\mathbf{x},t)$, initial–boundary conditions, and partial solution data $v_3(\mathbf{x},t)$ at the final time -- along with complete solution fields $v_1(\mathbf{x},t)$ and $v_2(\mathbf{x},t)$, the Alikhanov--PINN framework is trained to recover the fractional order $\alpha$, diffusion coefficient $\lambda_1$, reaction coefficient $\lambda_2$, and the full solution of the NfPDE.

For the inverse problem, the training dataset is extended by incorporating 30 randomly sampled terminal time solutions into the original forward dataset. 
Table~\ref{tab:nn_params_TFRD} lists the optimal neural network configurations such as depth, width, learning rate, IBC enforcement and unknown parameters along with the training durations for all test cases. 
The stopping criteria are defined using the maximum epoch count $m_{\mathrm{stage}}$, and a tolerance threshold on the relative $\mathcal{E}_2$-error $\epsilon$. 

The influence of the mesh grading parameter $\gamma$, is explored for the time-fractional reaction-diffusion equation (TFRD) inverse problem.
Figure~\ref{fig:ex4_tfrd} shows the  evolution of the parameters $(\alpha, \lambda_1,\lambda_2)$ on different graded meshes $(\gamma= 1, \frac{2}{\alpha}, \frac{3-\alpha}{\alpha})$. 
It is evident that nonuniform meshes ($\gamma=2/\alpha, (3-\alpha)/\alpha$) yield parameter estimates that converge more precisely to the true values than the uniform mesh ($\gamma=1$), highlighting the superior accuracy of nonuniform discretization in inverse TFRD problems.
Lastly, Table~\ref{tab:param_estimation_r} compares the ground truth and inferred parameters for 2D TFRD problems. 
The results confirm that all parameters are accurately estimated and that both error metrics diminish when using nonuniform meshes in place of uniform mesh.

%%%%%%%%%%%%%%%%%%%%%%%%%%%%%%%%%%%%%%%%%%%%
\begin{table}[htbp]
\centering
\scriptsize
\caption{Optimal network hyperparameters for the inverse 2D‐TFRD problem.}
\label{tab:nn_params_TFRD}
\begin{tabular}{l c}
\toprule
Parameter                             & Value               \\
\midrule
$L$                                   & 1                   \\
$N_i$                                 & 30                  \\
Optimizer                             & Adam                \\
IBCs constraints                      & Hard                \\
Learning rate (Adam)                  & $1e-04$   \\
Learning rate (unknown parameters)    &
  \begin{tabular}[t]{@{}l@{}}
    $\alpha = 1e-04$\\
    $\lambda_1 = 1e-04$\\
    $\lambda_2 = 1e-04$
  \end{tabular}                       \\
$m_{\mathrm{stage}}$                  & 8000                \\
\bottomrule
\end{tabular}
\end{table}

%%%%%%%%%%%%%%%%%%%%%%%%%%%%%%%
\begin{figure}[htbp]
  \centering
  \begin{subfigure}[b]{0.32\textwidth}
    \includegraphics[width=\textwidth]{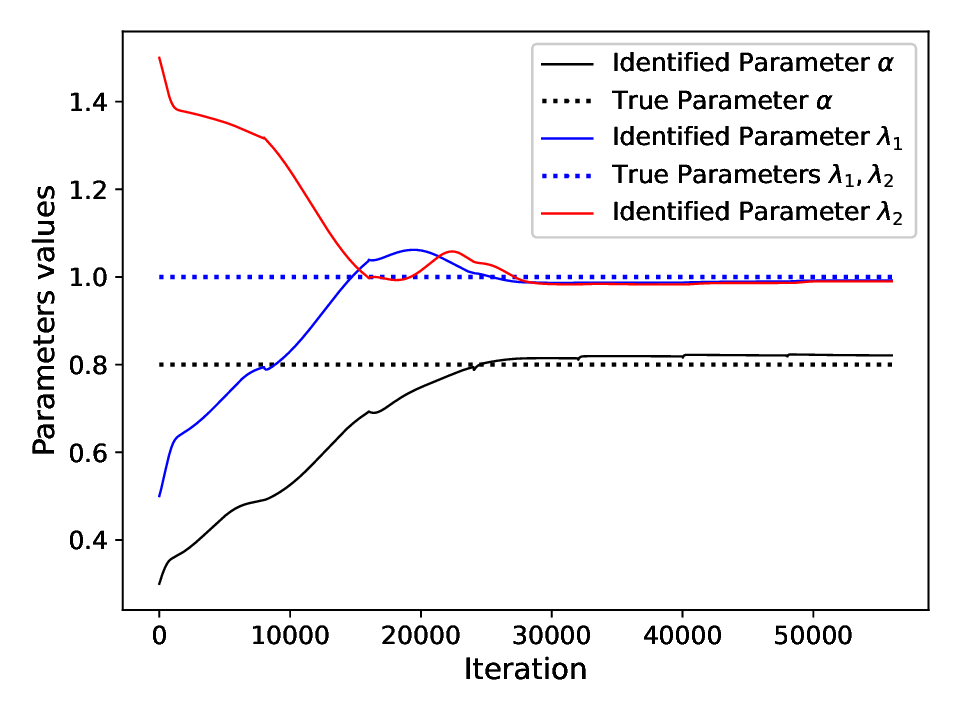}
    \caption{$\gamma=1$.}
    \label{fig:ex4_tfrd1}
  \end{subfigure}
  \hfill
  \begin{subfigure}[b]{0.32\textwidth}
    \includegraphics[width=\textwidth]{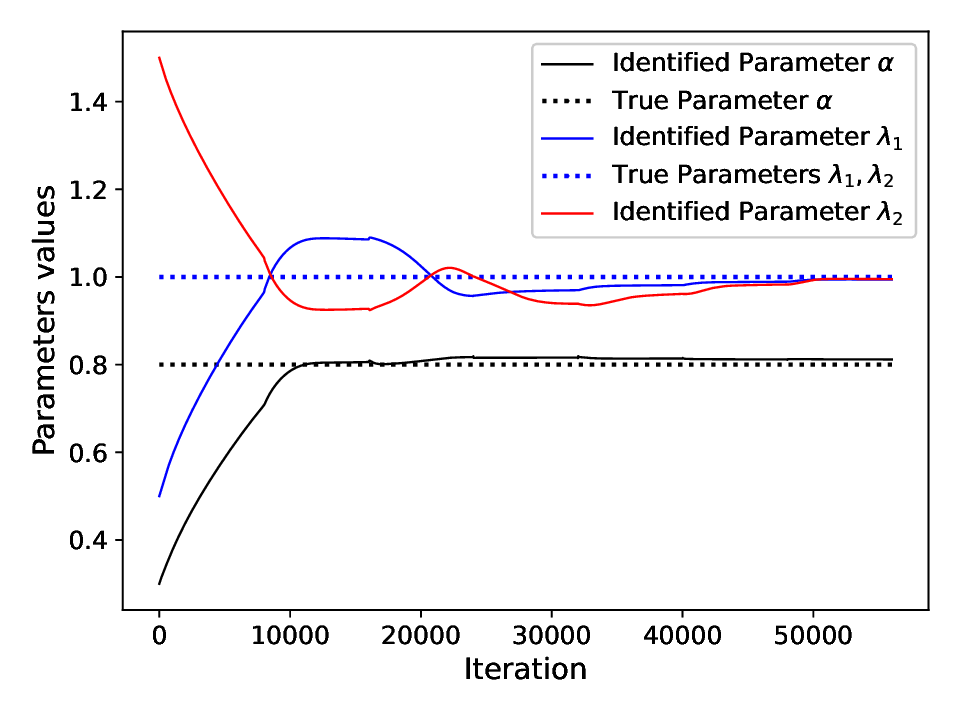}
    \caption{$\gamma=2/\alpha$.}
    \label{fig:ex4_tfrd2}
  \end{subfigure}
  \hfill
  \begin{subfigure}[b]{0.32\textwidth}
    \includegraphics[width=\textwidth]{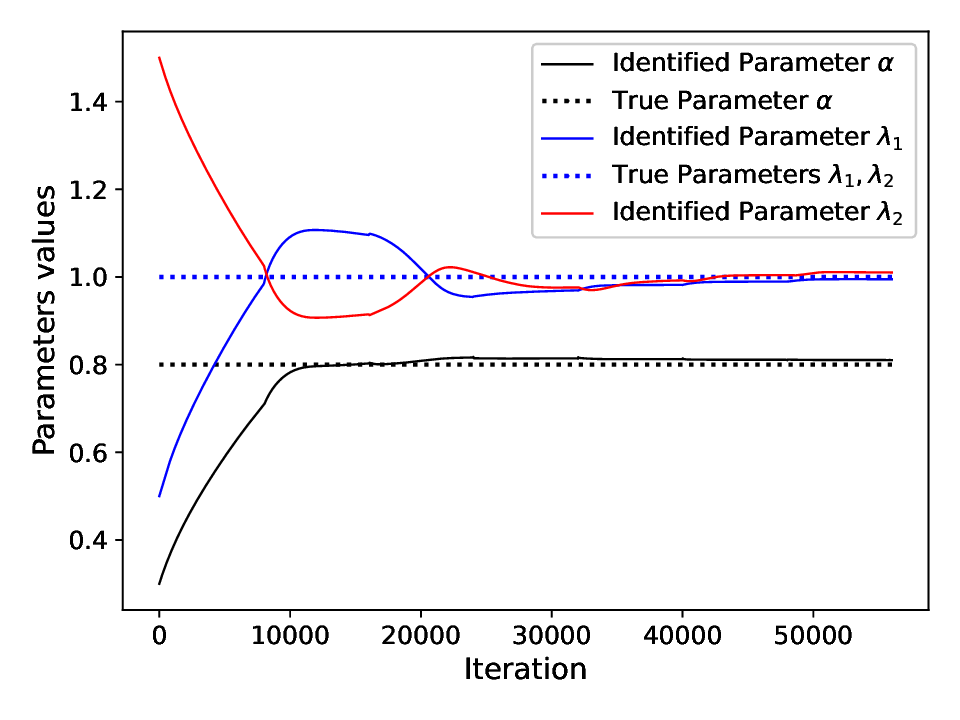}
    \caption{$\gamma=(3-\alpha)/\alpha$.}
    \label{fig:ex4_tfrd3}

  \end{subfigure}
  \caption{Evolution of the estimated parameters $\alpha$, $\lambda_1$, and $\lambda_2$ during training of the Alikhanov--XfPINNs for the 2D-TFRD equation on nonuniform meshes with varying grading parameter $\gamma$.}
  \label{fig:ex4_tfrd}
\end{figure}

%%%%%%%%%%%%%%%%%%%%%%%%%%%%%%%%%%
\begin{table}[h]
\centering
\small 
\caption{Comparison of true and estimated parameters, relative $L^2$-error, and maximum error for different $\gamma$‐values.}
\label{tab:param_estimation_r}
\begin{tabular}{c l l c c}
\hline
$\gamma$ & True Parameters & Estimated Parameters & $E_{2}(K)$ & $E_{\infty}(K)$ \\
\hline
$1$ &
\begin{tabular}[t]{@{}l@{}}
$\alpha = 0.8$\\
$\lambda_{1} = 1$\\
$\lambda_{2} = 1$
\end{tabular}
&
\begin{tabular}[t]{@{}l@{}}
$\alpha = 0.821$\\
$\lambda_{1} = 0.992$\\
$\lambda_{2} = 0.990$
\end{tabular}
& $8.299e-03$ & $6.951e-03$ \\
\hline
$2/\alpha$ &
\begin{tabular}[t]{@{}l@{}}
$\alpha = 0.8$\\
$\lambda_{1} = 1$\\
$\lambda_{2} = 1$
\end{tabular}
&
\begin{tabular}[t]{@{}l@{}}
$\alpha = 0.812$\\
$\lambda_{1} = 0.994$\\
$\lambda_{2} = 0.995$
\end{tabular}
& $5.182e-03$ & $3.083e-03$ \\
\hline
$(3-\alpha)/\alpha$ &
\begin{tabular}[t]{@{}l@{}}
$\alpha = 0.8$\\
$\lambda_{1} = 1$\\
$\lambda_{2} = 1$
\end{tabular}
&
\begin{tabular}[t]{@{}l@{}}
$\alpha = 0.810$\\
$\lambda_{1} = 0.995$\\
$\lambda_{2} = 1.010$
\end{tabular}
& $4.599e-03$ & $2.875e-03$ \\
\hline
\end{tabular}
\end{table}
\end{example}

%%%%%%%%%%%%%%%%%%%%%%%%%%%%%%%% EXAMPLE 3.5.
\begin{example}
Here, we investigate the 2D time-fractional Allen–Cahn (TFAC) equation for simultaneous parameter inference and solution reconstruction. 
The proposed methodology is applied on TFAC instances with a weak initial singularity to recover the unknown coefficients. 
Accuracy is quantified by the maximum-norm error $E_{\infty}(K)$ and the relative $L^2$-error $E_2(K)$. 
Throughout, $L$ denotes the network depth (number of layers), $N_i$ is the width of layer $i$.

Consider the inverse problem of TFAC equation on $\Omega^2=[0,L]\times [0,L]$
\begin{equation}\label{tfpf_inv}
\begin{aligned}
{} _{0}^{\mathcal{C}}\mathcal{\partial}_t^{\alpha}v(x,y,t)&=\Psi(\epsilon^2\Delta v(x,t)-f(v)) + g(x,y,t),\quad&(x,y,t)&\in\Omega^2\times[0,T],\\
v(x,y,0)&=\varphi(x,y),\quad&(x,y)&\in\Omega^2,\\
v(x,y,t)&=\phi(x,y,t),\quad&(x,y,t)&\in\partial\Omega^2\times[0,T],\\
v(x,y,T)&=\psi(x,y),\quad&(x,y)&\in\Omega^2,
\end{aligned}
\end{equation}
%\Tmatthias{what is the 2D domain $\Omega$ ? the notation $\Omega^2$ is not common, unless it is a cross-product of intervals.\\
%\textcolor{red}{Answer: The one-dimensional domain is $\Omega=[0,L]$ and corresponding one-dimensional domain is $\Omega^2=\Omega\times \Omega=[0,L]\times [0,L]$. }}
where, the parameters $\alpha$, $\epsilon$, and $\Gamma$ are treated as unknowns.
Given $f(v)$ reaction term, $g(x,y,t)$ source term, $\varphi(x,y)$ ICs, $\phi(x,y,t)$ BCs, and a sparse set of terminal-time observations $\psi(x,y)$ at $t=T$, Alikhanov--XfPINNs framework is employed to identify the fractional order $\alpha$, the interfacial width $\epsilon$, and the mobility parameter $\Psi$.

Table~\ref{tab:nn_params_TFAC} summarizes the network configurations used for the inverse TFAC problems. 
The listed parameters include the layer count $L$, the number of neurons $N_i$ per layer, the IBC enforcement type, the optimizer step size, and the target parameter set. 
The termination criteria are also defined by the stage-wise maximum iteration count $m_{\mathrm{stage}}$ and the convergence threshold $\varepsilon$.

We now evaluate the effect of time-mesh grading on the recovery of inverse parameters for the TFAC problem. 
For the two-dimensional TFAC equation with a weak initial singularity, Figure~\ref{fig:ex5_tfac} plots the Alikhanov--XfPINNs estimates of $\alpha$, $\epsilon$, and $\Gamma$ obtained with grading parameters $\gamma=1$, $\gamma=2/\alpha$, and $\gamma=(3-\alpha)/\alpha$. 
Compared to the uniform mesh ($\gamma = 1$), the graded settings $\gamma = 2/\alpha$ and $\gamma = (3-\alpha)/\alpha$ bring the estimates closer to the ground truth, underscoring the advantage of nonuniform time discretizations for the inverse TFAC problem. 
Table~\ref{tab13} presents and compares the exact and recovered parameters, reporting the maximum-norm error $E_{\infty}$ and the relative $L^2$-error $E_2(K)$. 
In line with the visual trends, both error measures are reduced on graded meshes, and all parameters are identified with high accuracy.

%%%%%%%%%%%%%%%%%%%%%%%%%%%%
\begin{table}[htbp]
\centering
\scriptsize
\caption{Optimal network hyperparameters for the inverse 2D‐TFAC problem.}
\label{tab:nn_params_TFAC}
\begin{tabular}{l c}
\toprule
Parameter                             & Value               \\
\midrule
$L$                                   & 1                   \\
$N_i$                                 & 30                  \\
IBCs constraints                      & Hard                \\
Optimizer                             & Adam                \\
Learning rate (Adam)                  & $1e-04$   \\
Learning rate (unknown parameters)    &
  \begin{tabular}[t]{@{}l@{}}
    $\alpha = 1e-04$\\
    $\varepsilon = 1e-04$\\
    $\Psi = 1e-04$
  \end{tabular}                       \\
$m_{\mathrm{stage}}$                  & 5000                \\
\bottomrule
\end{tabular}
\end{table}

%%%%%%%%%%%%%%%%%%%%%%%%%%%%
\begin{figure}[htbp]
  \centering
  \begin{subfigure}[b]{0.32\textwidth}
    \includegraphics[width=\textwidth]{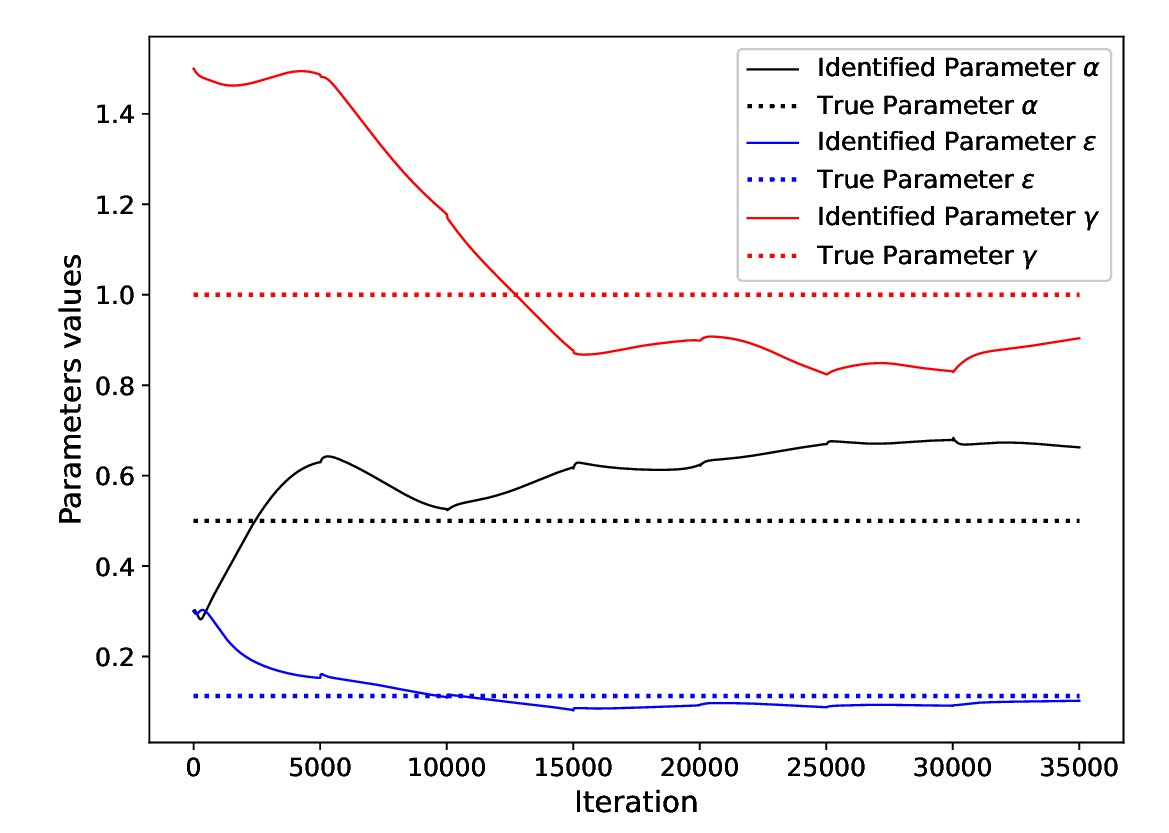}
    \caption{$\gamma=1$.}
    \label{fig:ex4_tfac1}
  \end{subfigure}
  \hfill
  \begin{subfigure}[b]{0.32\textwidth}
    \includegraphics[width=\textwidth]{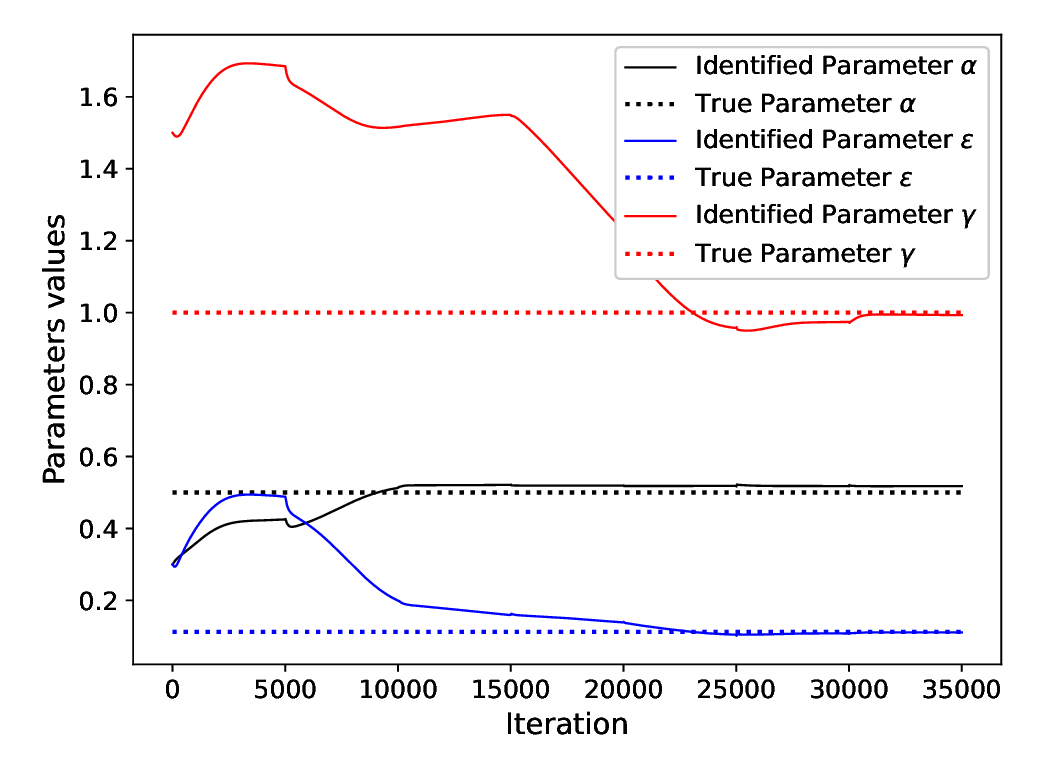}
    \caption{$\gamma=2/\alpha$.}
    \label{fig:ex4_tfac2}
  \end{subfigure}
  \hfill
  \begin{subfigure}[b]{0.32\textwidth}
    \includegraphics[width=\textwidth]{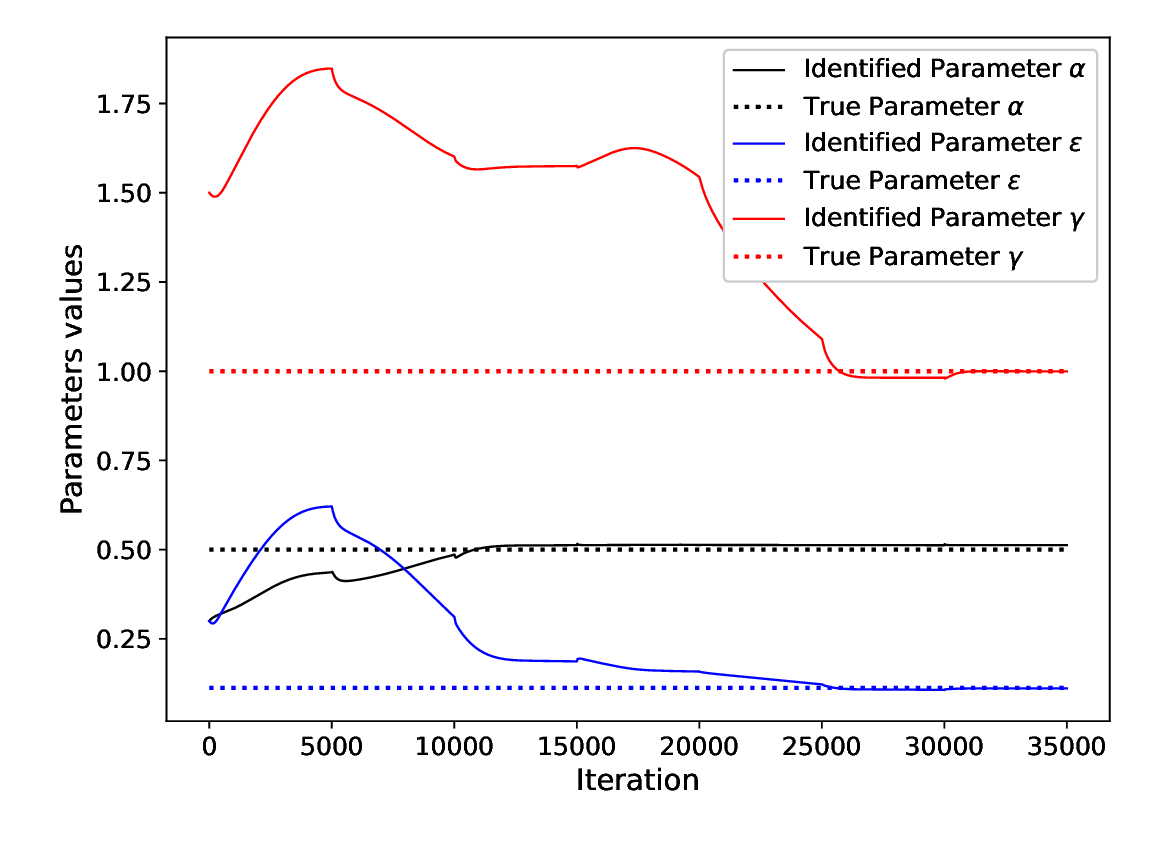}
    \caption{$\gamma=(3-\alpha)/\alpha$.}
    \label{fig:ex4_tfac3}

  \end{subfigure}
  \caption{Evolution of the estimated parameters $\alpha$, $\varepsilon$, and $\Gamma$ during training of the Alikhanov--XfPINNs for the 2D-TFAC equation on nonuniform meshes with varying grading parameter $\gamma$.}
  \label{fig:ex5_tfac}
\end{figure}

%%

%%%%%%%%%%%%%%%%%%%%%%%%%%%%%%%%%%update table
\begin{table}[htbp]
\centering
\small % Reduce font size to fit within page width
\caption{Comparison of true and estimated parameters, relative $L^2$-error, and maximum error for different $\gamma$‐values.}
\label{tab13}
\begin{tabular}{c l l c c}
\hline
$\gamma$ & True Parameters & Estimated Parameters & $E_{2}(K)$ & $E_{\infty}(K)$ \\
\hline
$1$ &
\begin{tabular}[t]{@{}l@{}}
$\alpha = 0.5$\\
$\varepsilon = \sqrt{2}/4\pi$\\
$\Psi = 1$
\end{tabular}
&
\begin{tabular}[t]{@{}l@{}}
$\alpha = 0.663$\\
$\varepsilon = 0.102$\\
$\Psi = 0.903$
\end{tabular}
& $9.572\text{e-}03$ & $1.258\text{e-}02$ \\
\hline
$2/\alpha$ &
\begin{tabular}[t]{@{}l@{}}
$\alpha = 0.5$\\
$\varepsilon = \sqrt{2}/4\pi$\\
$\Psi = 1$
\end{tabular}
&
\begin{tabular}[t]{@{}l@{}}
$\alpha = 0.518$\\
$\varepsilon = 0.111$\\
$\Psi = 0.993$
\end{tabular}
& $3.280\text{e-}03$ & $2.553\text{e-}03$ \\
\hline
$(3-\alpha)/\alpha$ &
\begin{tabular}[t]{@{}l@{}}
$\alpha = 0.5$\\
$\varepsilon = \sqrt{2}/4\pi$\\
$\Psi = 1$
\end{tabular}
&
\begin{tabular}[t]{@{}l@{}}
$\alpha = 0.513$\\
$\varepsilon = 0.111$\\
$\Psi = 0.999$
\end{tabular}
& $2.974\text{e-}03$ & $2.109\text{e-}03$ \\
\hline
\end{tabular}
\end{table}
\end{example}
%%%%%%%%%%%%%%%

%%%%%%%%%%%%%%%%%%%%%%%%%%%%%%%% Conclusion
\section{Conclusion}\label{sec6} % ME: rewritten
This work introduces a data-driven approach to solving nonlinear fractional partial differential equations (NfPDEs) that integrates an accelerated Alikhanov discretization on nonuniform time meshes with physics-informed neural networks (PINNs). 
This new method is called an \textit{Alikhanov-PINN}. 
This method builds efficient surrogate models that accurately approximate forward and inverse solutions of high-dimensional NfPDEs, including those with initial singularities and unknown analytical solutions. Unlike traditional finite difference-based fPINNs, our approach avoids full-domain training and dense system solvers. It employs a time-marching adaptive strategy that enables localized and efficient approximation on each temporal interval. 
Additionally, the method requires no spatial meshing, making it effective for irregular and curved domains in high dimensions. 
Including an adaptive activation function significantly accelerates convergence, and the flexible incorporation of hard and soft constraints enables the method to handle homogeneous and nonhomogeneous initial-boundary conditions.

We systematically examine the influence of network depth, adaptive activation, and collocation strategies on predictive accuracy for forward problems.
Extensive experiments demonstrate the robustness of the method across 1D, 2D, and geometrically complex domains and validate its scalability for long-time simulations. 
In the inverse problem setting, Alikhanov-PINNs reliably identify unknown parameters and accurately recover latent dynamics.
The auxiliary time-marching configuration, termed Alikhanov-fPINN--M, further enforces the discrete residual sequentially in time and provides an auditable protocol for examining empirical temporal convergence with respect to the maximal time step on graded meshes. 
Under controlled training tolerances, this configuration provides clear numerical evidence of the expected temporal convergence behavior of the underlying Alikhanov discretization.
The observed convergence behavior aligns with theoretical expectations, affirming the method’s numerical reliability. 
Looking ahead, this work establishes a foundation for several promising directions. 
These include extending the framework to accommodate other types of fractional operators, applying the method to broader classes of PDEs in physics and biology, and developing theoretical analyses of its approximation, generalization, and convergence properties. 
Furthermore, refining the methodology for hyperparameter selection and adaptive sampling strategies is essential for future optimization.

%\Tmatthias{just a comment, maybe for outlook: For future work, one could also use real data for some known model, and determine the best fitting order $\alpha$ (also as a function of time: $\alpha(t)$. This would (probably) provide a justification for using fractional time derivatives (some researchers have doubts about this).\\
%\textcolor{red}{Yes sir, I  agree that applying the proposed framework to real observational data and estimating the fractional order $\alpha$ directly from measurements would provide an important data-driven justification for the use of time-fractional models. 
%In particular, allowing the order to vary in time, i.e., $\alpha=\alpha(t)$, could help capture evolving memory effects and anomalous transport regimes that may not be well represented by a constant fractional order. Since the present work focuses on the numerical formulation and validation of the proposed Alikhanov-XfPINNs framework using benchmark problems, the incorporation of real data and the identification of constant or variable fractional orders are left as promising directions for future research  }}

In future work, we will apply the proposed framework to real observational data. Estimating the fractional order $\alpha$ directly from measurements would provide important data-driven justification for using time-fractional models. 
Allowing the order to vary over time, i.e., $\alpha=\alpha(t)$, could capture evolving memory effects and anomalous transport regimes that a constant fractional order may not adequately represent. 

% Since the present work focuses on the numerical formulation and validation of the proposed Alikhanov-XfPINNs framework using benchmark problems, the incorporation of real data and the identification of constant or variable fractional orders are left as promising directions for future research 

%%%%%%%%%%%%%%%%%%%%%%%%%%%%%%%%%%%%%%
%\section*{Author Contributions}
%Himanshu Kumar Dwivedi: Conceptualization, Soft\-ware, Methodology, Validation, Formal Analysis, Investigation, Writing-Original Draft, Writing-Review and Editing, Visualization.\\
%Matthias Ehrhardt: Conceptualization, Supervision, Methodology, Formal Analysis, Investigation, Writing-Review and Editing.\\
%Rajeev: Supervision, Funding acquisition, Writing-Review and Editing, Validation.

%%%%%%%%%%%%%%%%%%%%%%%%%%
%\section*{Data availability}
%The datasets supporting this study are available from the corresponding author upon reasonable request.

%%%%%%%%%%%%%%%%%%%%%%%%%%
%\section*{Declarations}
%\subsection*{Funding declarations}
%The authors extend their gratitude to the anonymous reviewers for their insightful comments and valuable  suggestions, which greatly enhanced the quality of this paper.
 %{The second author expresses gratitude for the project grant provided by SERB, Government of India, as outlined in Sanction Letter \\
 %MTR/2022/000169.}

%%%%%%%%%%%%%%%%%
%\subsection*{Conflict of interest}
%The authors confirm they have no conflicts of interest to report.

%%%%%%%%%%%%%%% Refs


\begin{thebibliography}{99}

\bibitem {P11} 
A. Krizhevsky, I. Sutskever, G. E. Hinton,
ImageNet classification with deep convolutional neural networks,
Adv. Neural Process. Syst. {(2012)} 1097-1105.

\bibitem {P12} 
B. M. Lake, R. Salakhutdinov, J. B. Tenenbaum, 
Human-level concept learning through probabilistic program induction, 
Science 350 {(2015)} 1332-1338.

\bibitem {P21} 
M. Raissi, P. Perdikaris, G. E. Karniadakis,  
Physics-informed neural networks: A deep learning framework for solving forward and inverse problem involving nonlinear partial differential equations, 
{ J. Comput. Phys.\/} {378} {(2019)} 686-707.

\bibitem {P22} 
G. E. Karniadakis, I. G. Kevrekidis, L. Lu, P. Perdikaris, S. Wang, L. Yang, 
Physics-informed machine learning,
{Nat. Rev. Phys. \/} {3(6)} {(2021)} 422-440.

\bibitem {P27} 
X. Jin, S. Cai, H. Li, G. E. Karniadakis, 
NSFnets (Navier-Stokes flow nets): 
Physics-informed neural networks for the incompressible Navier-Stokes equations, 
{ J. Comput. Phys.\/} {426} {(2021)} 109951.

\bibitem {P26} 
M. Raissi, Z. Wang, M. S. Triantafyllou, G. E.  Karniadakis, 
Deep learning of vortex-induced vibrations, {  J. Fluid Mech.\/} {861} {(2019)} 119-137.

\bibitem {P213} 
D. Liu, Y. Wang, 
Multi-fidelity physics-constrained neural network and its application in material modeling,
{  J. Mech. Des.\/} 141(12) {(2019)} 121403.

\bibitem {P24} 
{G. Kissas, Y. Yang, E. Hwuang, W. R. Witschey, J. A. Detre, P. Perdikaris}, 
Machine learning in cardiovascular flows modeling: predicting arterial blood pressure from non-invasive 4D flow MRI data using physics-informed neural networks, 
{  Comput. Methods Appl. Mech. Eng.\/} {358} {(2020)} 112623.

\bibitem {P215} 
{J. Sirignano, K. Spiliopoulos},
DGM: A deep learning algorithm for solving partial differential equations,
{ J. Comput. Phys. \/}  {375}  {(2018)}  1339-1364.

\bibitem {P219}
{L. Lu, X. Meng, Z. Mao, G. E. Karniadakis}, 
DeepXDE: a deep learning library for solving differential equations, {SIAM Rev. \/}  {63(1)}  {(2021)} 208-228.

\bibitem {P220} 
{ S. Wang, S. Sankaran, P. Perdikaris},
Respecting causality for training physics-informed neural networks, 
{Comput. Methods Appl. Mech. Eng Phys. \/}  {421}  {(2024)} 116813.


\bibitem {P113}
{ C. Basdevant, M. Deville, P. Haldenwang, J. Lacorix, J. Ouazzani, R. Peyret, P. Orlandi, A. Patera},
Spectral and finite difference solutions of the Burgers' equation,
{ Comput. Fluids. \/}  {14}  {(1986)} 23-41.

\bibitem {P410}
{Z. Sun, G. Gao}, 
Fractional differential equations: finite difference methods, 
{ De Gruyter, Berlin \/}  {(2020)}.

\bibitem {jin} 
{  B. Jin, R. Lazarov,  Z.  Zhou}, 
An analysis of the L1 scheme for the subdiffusion equation with nonsmooth data, 
{ IMA J. Numer. Anal. \/}  {36}  {(2016)}  197–221.

 \bibitem {stynes2} 
 { M. Stynes, E. O´Riordan,  J. Gracia}, 
 Error analysis of a finite difference method on graded meshes for a time-fractional diffusion equation,
 { SIAM J. Numer. Anal. \/}  { 55(2) }  {(2017)} 1057–1079.

 \bibitem {lean} 
 { W. McLean}, 
 Regularity of solutions to a time-fractional diffusion equation,
 { ANZIAM J. \/} {52}  {(2010)} 123–138.

  
\bibitem {stynes1} 
{H. Chen, M. Stynes}, 
Error analysis of a second-order method on fitted meshes for a time-fractional diffusion problem, 
{J. Sci. Comput. \/}  {79(1)}  {(2019)}  624–647. 

\bibitem {sjiang2017}  
{S. Jiang, J. Zhang, Q. Zhang, Z. Zhang},   
Fast evaluation of the Caputo fractional derivative and its applications to fractional differential equations, 
{Commun. Comput. Phys.\/} {21(3)} {(2017)}. 650-678.

\bibitem {dwivedi1}  
{H. K. Dwivedi, Rajeev},
A novel fast second order approach with high-order compact difference scheme and its analysis for the tempered fractional Burgers equation, 
{Math. Comput. Simul.\/} {227} {(2025)} 168-188. 

\bibitem {dwivedi_2} 
{H. K. Dwivedi, Rajeev}, 
A novel fast tempered algorithm with high-accuracy scheme for 2D tempered fractional reaction-advection-subdiffusion equation, 
{Comput. Math. Appl. \/} {176}  {(2024)} 371-397.

 \bibitem {fast_tempered}  
 J. Cao, A. Xiao, W. Bu,  
 Finite difference/finite element method for the tempered time fractional advection-dispersion equation with fast evaluation of Caputo derivatives,
 { J. Sci. Comput. \/} {}{(2020)} 83:48.

\bibitem {secondfast}  
{Y. Yan, Z. Z. Sun, J. Zhang},
Fast evaluation of the Caputo fractional derivative and its applications to fractional diffusion equations: a second order scheme, 
{ Commun. Comput. Phys. \/} {22(4)} {(2017)} 1028-1048.

\bibitem {fast_non} 
{X. Li, H. L. Liao, L. Zhang}, 
A second-order fast compact scheme with unequal time-steps for subdiffusion problems, 
{Numer. Algor.  \/} {86} {(2021)} 1011–1039.

\bibitem{P417} 
Y. Luchko, W. Rundell, M. Yamamoto, L. Zuo, 
Uniqueness and reconstruction of an unknown semilinear term in a time-fractional reaction diffusion equation,
{Inverse Probl. \/}  {29(6)}  {(2013)}  065019.

\bibitem {P418}
{M. Kern}, 
Numerical methods for inverse problems,
{Wiley, Hoboken \/} {(2016)}.

\bibitem {P419} 
D. Lukyanenko, R. Argun, A. Borzunov, A. Gorbachev, V. Shinkarev, M. Shislenin, A. Yagola, 
On the features of numerical solution of coefficient inverse problems for nonlinear equations of the reaction-diffusion-advection type with data of various types, 
{ Differ. Equ. \/}  {59(12)}  {(2023)}  1734-1757

\bibitem {P420} 
{ X. An, Q. Wang, F. Liu, V. V. Anh, I. W. Turner}, 
Parameter estimation for time-fractional Black-Scholes equation with S\&P 500 index option, { Numer. Algor. \/}  {95(1)}  {(2024)} 1-30.

\bibitem {P425} 
{ M. Srati, A. Oulmelk, L. Afraites, A. Hadri, M. A. Zaky, A. Aldraiweesh, A. S. Hendy}, 
An inverse problem of determining the parameters in diffusion equations by using fractional physics-informed neural networks,
{ Appl. Numer. Math. \/}  {208}  {(2025)} 189-213.

\bibitem {P427} 
L. Ma, R. Li, F. Zeng, L. Guo, G. E. Karniadakis,
Bi-orthogonal fPINN: a physics-informed neural network method for solving time-dependent stochastic fractional PDEs, 
{ Commun. Comput. Phys. \/}  {34(4)}  {(2023)} 1133-1176.

\bibitem {P429} 
X. Fang, L. Qiao, F. Zhang, F. Sun, 
Explore deep network for a class of fractional partial differential equations, 
{ Chaos Solitons Fractals \/}  {172}  {(2022)} 113528.

\bibitem {P4}
{ J. Shi, X. Liu, X, Yang}, 
Data-driven solutions and parameter estimation of the high-dimensional time-fractional reaction-diffusion equations using an improved fPINN method,
{ Nonlinear Dyn. \/}  {113}  {(2025)} 9577-9604.

\bibitem {P431} 
{G. Pang, L. Lu, G. E. Karniadakis}, 
fPINNS: fractional physics-informed neural networks,
{ SIAM J. Sci. Comput. \/}  {41(4)}  {(2019)} 2603-2626 .

\bibitem {P433}
{ S. Wang, H. Zhang, X. Jiang}, 
Fractional physics-informed neural networks for time-fractional phase field models,
{ Nonlin. Dyn. \/}  {110(3)}  {(2022)} 2715-2739.

\bibitem {P435} 
{ Z. Ma, J. Hou, W. Zhu, Y. Peng, Y. Li},
PMNN: physical model driven neural network for solving time-fractional differential equations, 
{ Chaos Solitons Fractals \/}  {177}  {(2023)} 114238.

\bibitem {P434} 
H. Ren, X. Meng, R. Liu, J. Hou, Y. Yu,
A class of improved fractional physics-informed neural networks,
{ Neurocomputing \/}  {562}  {(2023)} 126890.

\bibitem {liao16}
{ H. L. Liao, , D. Li,  J. Zhang}, 
Sharp error estimate of nonuniform L1 formula for time-fractional reaction-subdiffusion equations, 
{ SIAM J. Numer. Anal. \/} {56} {(2018)} 1112–1133.

\bibitem {Alikhanov}  
{ A. A. Alikhanov},
A new difference scheme for the time fractional diffusion equation, 
{ J. Comput. Phys.  \/} {280}, {(2015)} 424-438.

\bibitem {Liao_arxiv} 
{ H. L. Liao,  W. Mclean, J. Zhang},
A second-order scheme with nonuniform time steps for a linear reaction-subdiffusion problem,
{Commun. Comput. Phys. \/}  {30(2)}, {(2021)} 567-601.

\bibitem{jagtap} 
{A. D. Jagtap, K. Kawaguchi, G. E. Karniadakis},
Adaptive activation functions accelerate convergence in deep and physics-informed neural networks,
{J. Comput. Phys. \/}  {404}, {(2020)} 109136.

\bibitem {geonet}
{H. Gao, L. Sun, J. Wang}, 
PhyGeoNet: physics-informed geometry-adaptive convolutional neural networks for solving parametrized steady-state PDE on irregular domain,
{J. Comput. Phys. \/}  {428}, {(2021)} 110079.

\bibitem{shi} 
{J. Shi, X. Liu, X. Yang}, Data-driven solutions and parameter estimation of the high-dimensional time-fractional reaction-diffusion equations using an improved fPINN method,
{Nonlin. Dyn. \/}  {113}  {(2025)}9577-9604.

\bibitem {tensorflow} 
{M. Abadi, A. Agarwal, P. Barham, E. Brevdeo, Z. Chen, C. Citro, G. S. Corrado, A. Davis, J. Dean, M. Devin, et al.},
Tensorflow: Large-scale machine learning on heterogeneous distributed systems,
{arXiv:1603.04467 \/} {(2016)}.

\bibitem {adam} 
{A. Barkat, P. Bianchi},
Convergence and dynamical behavior of the ADAM algorithm for non-convex stochastic optimization,
{SIAM J. Optim. \/} {31(1)}  {(2021)} 244-274.

\end{thebibliography}
\end{document}